\documentclass[10pt]{article}
\usepackage{bbm}
\usepackage{amsfonts}
\usepackage{amsmath}
\allowdisplaybreaks
\usepackage{amssymb}
\usepackage{fancyhdr}
\usepackage{setspace}
\usepackage{latexsym}
\usepackage{mathrsfs}
\usepackage{amsthm}

\usepackage{cite,enumitem,graphicx}

\usepackage[colorlinks=true,urlcolor=blue,
citecolor=red,linkcolor=blue,linktocpage,pdfpagelabels,
bookmarksnumbered,bookmarksopen]{hyperref}
\usepackage[english]{babel}

\numberwithin{equation}{section}
\newtheorem{theorem}{Theorem}[section]
\newtheorem{proposition}[theorem]{Proposition}
\newtheorem{lemma}[theorem]{Lemma}

\newtheorem{remark}[theorem]{Remark}

\newcommand{\ud}{\mathrm{d}}

\makeatletter 
\renewcommand\theequation{\thesection.\arabic{equation}}
\@addtoreset{equation}{section}
\makeatother
\numberwithin{equation}{section}

\begin{document}

\begin{center}
\textbf{Normalized solutions to mixed dispersion nonlinear Schr\"odinger system with coupled nonlinearity}\\
\end{center}

\begin{center}
Zhen-Feng Jin$^{1}$, Guotao Wang$^{1}$, Weimin Zhang$^{2,*}$\\

\smallskip
$^{1}$ Shanxi Key Laboratory of Cryptography and Data Security \& School of Mathematical Sciences, Shanxi Normal University, Taiyuan, Shanxi, 030031, P.R. China\\

\smallskip
$^{2}$ School of Mathematical Sciences, Zhejiang Normal University, Jinhua, Zhejiang, 321004, P.R. China
\end{center}

\begin{center}
\renewcommand{\theequation}{\arabic{section}.\arabic{equation}}
\numberwithin{equation}{section}
\footnote[0]{\hspace*{-7.4mm}
$^*$Corresponding author.
{E-mail addresses: jinzhenfeng@sxnu.edu.cn (Z.~Jin),
wgt2512@163.com (G. Wang), zhangweimin2021@gmail.com (W.~Zhang).}}
\end{center}

\begin{abstract}
In this paper, we consider the existence of normalized solutions for the following biharmonic nonlinear Schr\"{o}dinger system
\begin{equation*}
\begin{cases}
\Delta^2u+\alpha_{1}\Delta u+\lambda u=\beta r_{1}|u|^{r_{1}-2}|v|^{r_{2}} u & \text { in } \mathbb{R}^{N}, \\
\Delta^2v+\alpha_{2}\Delta v+\lambda v=\beta r_{2}|u|^{r_{1}}|v|^{r_{2}-2} v & \text { in } \mathbb{R}^{N}, \\
\int_{\mathbb{R}^{N}} (u^{2}+v^{2})\ud x=\rho^{2},
\end{cases}
\end{equation*}
where $\Delta^2u=\Delta(\Delta u)$ is the biharmonic operator, $\alpha_{1}$, $\alpha_{2}$, $\beta>0$, $r_{1}$, $r_{2}>1$, $N\geq 1$. $\rho^2$ stands for the prescribed mass, and $\lambda\in\mathbb{R}$ arises as a Lagrange multiplier. Such single constraint permits mass transformation in two materials.
When $r_{1}+r_{2}\le 2+\frac{8}{N}$, we obtain a dichotomy result  with respect to the mass for the existence of nontrivial ground states. Especially when $\alpha_1=\alpha_2$, the ground state exists for all $\rho>0$ if and only if $r_1+r_2<\min\left\{\max\left\{4, 2+\frac{8}{N+1}\right\}, 2+\frac{8}{N}\right\}$. When $r_{1}+r_{2}\in\left(2+\frac{8}{N}, \frac{2N}{(N-4)^{+}}\right)$ and $N\geq 2$, we obtain the existence of radial nontrivial mountain pass solution for sufficiently small $\rho>0$.
\end{abstract}


\section{Introduction}
This paper is concerned with the existence of solutions for the biharmonic nonlinear Schr\"{o}dinger system
\begin{equation}\label{23082801}
\begin{cases}
\Delta^2u+\alpha_{1}\Delta u+\lambda u=\beta r_{1}|u|^{r_{1}-2}|v|^{r_{2}} u & \text { in } \mathbb{R}^{N}, \\
\Delta^2v+\alpha_{2}\Delta v+\lambda v=\beta r_{2}|u|^{r_{1}}|v|^{r_{2}-2} v & \text { in } \mathbb{R}^{N},
\end{cases}
\end{equation}
under the mass constraint
\begin{equation}\label{23101001}
\begin{aligned}
\int_{\mathbb{R}^{N}} (u^{2}+v^{2})\ud x=\rho^{2},
\end{aligned}
\end{equation}
where $\alpha_{1}$, $\alpha_{2}$, $\beta$, $\rho>0$, $N\geq 1$, $r_{1}$, $r_{2}>1$, $r:=r_{1}+r_{2}<2^{**}$, $\Delta^2u=\Delta(\Delta u)$ is the biharmonic operator, and $\lambda\in\mathbb{R}$ arises as a Lagrange multiplier, which is unknown. Here
$$2^{**}:=\frac{2N}{(N-4)^+},\, \mbox{namely $2^{**}=\frac{2N}{N-4}$ if $N\geq 5$, and $2^{**}=+\infty$ if $N\leq 4$},$$
which is called Sobolev critical exponent when $N\ge 5$. Any solution $(u, v, \lambda)$ of \eqref{23082801} satisfying \eqref{23101001} is usually called a normalized solution of \eqref{23082801}.

\smallskip
In recent years, many researchers have considered the normalized solutions for the following biharmonic Schr\"odinger equation
\begin{equation}\label{24101405}
\begin{aligned}
\begin{cases}
\Delta^2u+\alpha\Delta u+\lambda u=f(u) & \text { in } \mathbb{R}^{N},\\
\int_{\mathbb{R}^{N}} u^{2}\ud x=\rho^{2},
\end{cases}
\end{aligned}
\end{equation}
where $\alpha\in\mathbb{R}$, $\rho>0$ and $\lambda$ is a Lagrange multiplier. There is a fruitful research background on the biharmonic nonlinear Schr\"odinger equation. On the one hand, in order to regularize and stabilize the solutions to the Schr\"{o}dinger equation
\[
i\partial_t\psi+ \Delta\psi + f(\psi)=0 \quad \mbox{ in }\mathbb{R}^N\times (0,\infty),
\]
Karpman and Shagalov (see \cite{KS2000} and references therein) proposed the fourth order equation
\begin{equation}\label{2501182234}
i\partial_t\psi-\gamma \Delta^2\psi +\beta \Delta\psi + f(\psi)=0 \quad \mbox{ in }\mathbb{R}^N\times (0,\infty).
\end{equation}
On the other hand, to prevent blow-up in finite time, Fibich, Ilan and Papanicolaou \cite{FI_SIAM2002} added a small fourth-order dispersion
term in nonlinear Schr\"odinger equation as a nonparaxial correction.  We omit the details about these introductions, one can refer to \cite{Bonheure18, Bonheure192}.

\smallskip
Here we are concerned with the standing wave solutions, namely $\psi$ of the form
\begin{equation}\label{2501182233}
\psi(x, t)=e^{i\lambda t} u(x).
\end{equation}
If $f(\psi)=e^{i\lambda t}f(u)$, up to a scaling, \eqref{2501182234}-\eqref{2501182233} with prescribed $L^2$-norm can be reduced to the equation \eqref{24101405}. Note that the biharmonic term and Laplacian term are two dispersive terms, and a lack of homogeneity occurs provided $\alpha\neq 0$, which brings the main difficulties in deriving the existence of solutions. An approach to obtain solutions of \eqref{24101405} is variational method. That is, the critical points of the associated energy functional
\[
E(u)=\frac{1}{2}\|\Delta u\|_2^2-\frac{\alpha}{2}\|\nabla u\|_2^2-\int_{\mathbb{R}^N}\int_{0}^uf(t)\ud t \ud x, \quad \forall u\in \Lambda_{\rho},
\]
correspond to the solutions of \eqref{24101405}, where the constraint
\[
\Lambda_{\rho}:=\{u\in H^2(\mathbb{R}^N): \|u\|_2^2=\rho^2\},
\]
and the Sobolev space $H^{2}(\mathbb{R}^{N})$ is defined as follows
\begin{equation*}
\begin{aligned}
H^{2}(\mathbb{R}^{N}):=\left\{u\in L^{2}(\mathbb{R}^{N}):\nabla u,\Delta u\in L^{2}(\mathbb{R}^{N})\right\}
\end{aligned}
\end{equation*}
endowed with the equivalent norm $\|u\|:=\left(\|\Delta u\|^2_2+\|u\|^2_2\right)^{\frac{1}{2}}$. Observe that the biharmonic term $\|\Delta u\|_2^2$ possesses dominating role in energy functional when $\alpha\le 0$, while the Laplacian term has great effect on the geometry of energy functional when $\alpha>0$, especially which shows a convex-concave shape when we rescale the function $u$. Regardless of $\alpha\le 0$ or $\alpha>0$, there exist a lot of works concentrated on the nonlinearities $f(u)=\mu|u|^{q-2}u+|u|^{p-2}u$ with $\mu\geq0$.

\smallskip
The simplest case is the power nonlinearities $f(u)=|u|^{p-2}u$. Consequently the energy functional is
\[
\widetilde{E}(u)=\frac{1}{2}\|\Delta u\|_2^2-\frac{\alpha}{2}\|\nabla u\|_2^2-\frac1p\int_{\mathbb{R}^N}|u|^p \ud x, \quad \forall u\in \Lambda_{\rho}.
\]
It is convenient to introduce the parameter $\bar{r}:=2+\frac{8}{N}$, which is usually called  $L^2$-critical exponent or mass critical exponent for fourth order equation. By Gagliardo-Nirenberg inequality (see Lemma \ref{L23083006} below) and the inequality \eqref{23092305}, or by dilations, the sign of $p-\bar{r}$ decides the geometry of the energy functional $\widetilde{E}|_{\Lambda_{\rho}}$. Observe that $\widetilde{E}$ is coercive on $\Lambda_{\rho}$ if $2<p<\bar{r}$; and $\widetilde{E}$ is unbounded from below on $\Lambda_{\rho}$ if $\bar{r}<p\le 2^{**}$.

\smallskip
Many researchers treated the problem \eqref{24101405} with $f(u)=|u|^{p-2}u$ according to the value of $p$ and the sign of $\alpha$.
\begin{itemize}
\item  When $\alpha\leq 0$, Bonheure {\it et al.} \cite{Bonheure18} considered the coercive case, namely $p\in(2,\bar{r})$. Observe that the minimization level
\[
\overline{m}(\rho):=\underset{u\in \Lambda_\rho}{\inf} \widetilde{E}(u)
\]
is obviously sub-additive since the energy is invariant under translation. Whenever $\overline{m}(\rho)$ is negative, the vanishing of the minimizing sequence will not occur, hence $\overline{m}(\rho)$ can be attained. They proved that $\overline{m}(\rho)<0$ if $2<p<2+\frac{4}{N}$, so $\overline{m}(\rho)$ is attainable for all $\rho>0$. For $2+\frac{4}{N}\le p<\bar{r}$, $\overline{m}(\rho)$ could be zero if the mass is large. As a result, they obtained a dichotomy result with respect to mass.
\item Bonheure {\it et al.} \cite{Bonheure192} investigated the case $\alpha<0$ and $p\in[\bar{r}, 2^{**})$. They proved the existence of the ground states when $\rho\in (c_1,c_2)$ for two numbers $c_1\ge 0,\, c_2\in (c_1, \infty]$ depending only on $N,\, p$. In particular, there holds that $c_1=0$ if $p>\bar{r}$. Moreover, they constructed minimax levels by $\mathbb{Z}_2$-genus, obtaining the multiplicity of radial normalized solutions.
\item When $\alpha>0$, the problem is more involved since the Laplacian has an effect on the shape of energy functional.
Fern\'{a}ndez {\it et al.} \cite{Fernandez22} established some non-homogeneous Gagliardo-Nirenberg inequalities, and showed some results about the existence and non-existence of global minimizers for $p\in(2,\bar{r}]$ and local minimizers for $p\in(\bar{r},2^{**})$. Furthermore, when $p\in(\bar{r},2^{**})$, Luo and Yang \cite{Luo23} proved that \eqref{24101405} has another mountain pass type solution as $\rho$ is small.
\end{itemize}
   Next we introduce some recent results about $f(u)=\mu|u|^{q-2}u+|u|^{p-2}u$ with $\mu>0$.
\begin{itemize}
\item When $\alpha=0$, $N\ge 5$ and $p=2^{**}$, Ma and Chang \cite{Ma23} showed that for $q\in(2,\bar{r})$, \eqref{24101405} has a ground state solution provided $\rho$ is small; Liu and Zhang \cite{LiuJ23} proved that for $q\in(2,\bar{r})$, there exists a mountain pass type solution of \eqref{24101405} for small $\rho$.
\item When $\alpha=0$, $N\ge 5$, $\bar{r}\le q<p\leq2^{**}$, Chang {\it et al.} \cite{Chang23} discussed the existence, non-existence of normalized solutions for \eqref{24101405}, and they also proved the strong instability of standing waves.
\end{itemize}
Finally, for more general $L^2$-subcritical nonlinearities and $\alpha\in\mathbb{R}$, Luo and Zhang \cite{LuoH22} proved a dichotomy result with respect to the mass for the existence of ground states to \eqref{24101405}. Chen and Chen \cite{ChenJ23} considered nonlinearities $f$ involving Hardy-Littlewood-Sobolev upper critical and combined nonlinearities.

\smallskip
As far as we are aware, it seems that there is no literature considering normalized solutions involving mixed dispersion nonlinear Schr\"odinger system. Motivated by the above works, we are interested in the existence of normalized solutions for problem \eqref{23082801}-\eqref{23101001}. That is, we search for $(u, v, \lambda)$ satisfying \eqref{23082801}-\eqref{23101001}. A solution $(u, v, \lambda)$ is said {\it nontrivial}, which means that  $u\neq 0$ and $v\neq 0$. The functional
\begin{equation*}
\begin{aligned}
I(u,v)=\frac{1}{2}\int_{\mathbb{R}^{N}}(|\Delta u|^{2}+|\Delta v|^{2}) \ud x
-\frac{\alpha_{1}}{2}\int_{\mathbb{R}^{N}}|\nabla u|^{2} \ud x
-\frac{\alpha_{2}}{2}\int_{\mathbb{R}^{N}}|\nabla v|^{2} \ud x
-\beta\int_{\mathbb{R}^{N}}|u|^{r_{1}}|v|^{r_{2}} \ud x
\end{aligned}
\end{equation*}
is the corresponding variational functional of problem \eqref{23082801}-\eqref{23101001} defined on the constraint
\begin{equation}\label{2501181716}
\begin{aligned}
S_{\rho}:=\left\{(u,v)\in H^2(\mathbb{R}^{N})\times H^2(\mathbb{R}^{N})
:\|u\|^2_{2}+\|v\|^2_{2}=\rho^2\right\}.
\end{aligned}
\end{equation}
To the best of our knowledge, very few works concerned the single constraint for $u$ and $v$ in nonlinear Schr\"odinger system. However the constraint \eqref{2501181716} could be often encountered in physical world. In fact, $(u, v)\in S_\rho$ represents that the total mass of a system is conserved, but the mass of $u$ and $v$ may transform mutually, where the transformation may be aroused due to chemical reactions or the movement in physical space.

In this paper, we are focused on $\alpha_1, \alpha_2, \beta>0$. Note that the energy functional $I$ under the different cases $r\in (2, \bar{r})$, $r=\bar{r}$ and $r\in (\bar{r}, 2^{**})$ has different geometry on $S_\rho$. If either $r\in(2,\bar{r})$ or $r=\bar{r}$ and $\rho<\left(\frac{1}{2\mathcal{D}_1\beta}\right)^{\frac{N}{8}}$ where
\begin{equation}\label{2507061525}
\mathcal{D}_{1}=\left(\frac{r_1}{r}\right)^{r_1}\left(\frac{r_2}{r}\right)^{r_2} {C}_{N,r}^r,
\end{equation}
and ${C}_{N,r}$ is the optimal constant for Gagliardo-Nirenberg inequality \eqref{2411232346}, we get in Lemma \ref{l24070902} that $I$ is coercive on $S_\rho$. Hence we ask whether the minimization problem
\begin{equation}\label{2411071036}
\begin{aligned}
m(\rho):=\inf_{(u,v) \in S_\rho} I(u,v)
\end{aligned}
\end{equation}
 is attainable. Usually, normalized solutions with least energy are called normalized {\it ground states} (for short, we omit the word ``normalized"), so any minimizer of \eqref{2411071036} is a  ground state. As the method of \cite{Fernandez22, LuoH22}, the key ingredient is to establish the compactness of minimizing sequence. We emphasize here that $\alpha_{1}$, $\alpha_{2}>0$ and the fact that one can not obtain the strong convergence $\nabla u_{n}\rightarrow\nabla u$ in $L^{2}(\mathbb{R}^{N})$ only from $u_{n} \rightharpoonup u$ in $H^{2}(\mathbb{R}^{N})$ brings the main difficulties for problem \eqref{23082801}-\eqref{23101001}. For that, we first establish a strict sub-additivity for $m(\rho)$ in Lemma \ref{L24100202}. With this result and Lions' concentration compactness principle, we provide in Theorem \ref{T24100203} an alternative result for minimizing sequence, that is, the minimizing sequence either vanishes or has a convergent subsequence.

 \smallskip
 On the other hand, we consider another auxiliary minimization problem
\begin{equation}\label{2411111354}
m^{J}(\rho):=\underset{(u, v)\in S_{\rho}}{\inf} J(u, v),
\end{equation}
where
\begin{equation*}
J(u, v)=\frac{1}{2}\left(\|\Delta u\|^{2}_{2}+\|\Delta v\|^{2}_{2}\right)
-\frac{\alpha_{1}}{2}\|\nabla u\|^{2}_2-\frac{\alpha_{2}}2
\|\nabla v\|^{2}_2,\quad \forall (u, v)\in S_\rho.
\end{equation*}
It can be seen that $J$ is coercive on $S_\rho$, so $m^{J}(\rho)$ is always well-defined. We will show in Lemma \ref{24100101} that $m^{J}(\rho)$ is never achieved. If the minimizing sequence for $m(\rho)$ vanishes, there must hold $m(\rho)=m^J(\rho)$. Therefore we can rule out the occurrence of vanishing  if $m(\rho)<m^J(\rho)$. So the key
step is to establish the comparison between $m(\rho)$ and $m^J(\rho)$. As the ideas in  \cite{Fernandez22}, this comparison can be reduced into whether the supremum of
\begin{equation}\label{2411192136}
Q(u, v):= \frac{\int_{\mathbb{R}^{N}}|u|^{r_{1}}|v|^{r_{2}} \ud x}{\left(\left\|\left(\Delta + \frac{\alpha_1}{2}\right)u\right\|_2^2+ \left\|\left(\Delta + \frac{\alpha_2}{2}\right)v\right\|_2^2+\frac{\alpha_1^2-\alpha_2^2}{4}\|v\|_2^2\right)\left(\|u\|_2^2+\|v\|_2^2\right)^{\frac{r}{2}-1}}
\end{equation}
in $H^2(\mathbb{R}^N)\times H^2(\mathbb{R}^N)$ is finite, see subsection \ref{2501182319}. After careful analysis, we obtain the following dichotomy result.
\begin{theorem}\label{T24070901}
Assume that $\alpha_{1}$, $\alpha_{2}$, $\beta>0$, $r_{1}$, $r_{2}>1$.
Let $\bar{r}=2+\frac{8}{N}$, $r=r_1+r_2$ and $\mathcal{D}_{1}$ be given in \eqref{2507061525}. Then there exists some
\[
\begin{aligned}
\rho^*\in
\begin{cases}
 [0, \infty)\quad &\mbox{if}\;\; r<\bar{r},\\
\left[0, \left(\frac{1}{2\mathcal{D}_1\beta}\right)^{\frac{N}{8}}\right)\quad &\mbox{if}\;\; r=\bar{r},
\end{cases}
\end{aligned}
\]
such that $m(\rho)$ can be attained  if
\[
\begin{aligned}
\rho\in
\begin{cases}
 (\rho^*, \infty)\quad &\mbox{if}\;\; r<\bar{r},\\
\left(\rho^*, \left(\frac{1}{2\mathcal{D}_1\beta}\right)^{\frac{N}{8}}\right)\quad &\mbox{if}\;\; r=\bar{r},
\end{cases}
\end{aligned}
\]
where the minimizers are nontrivial. Moreover, \eqref{2411071036} can never be attained  if $\rho<\rho*$.
\end{theorem}

Notice that it is possible that $\rho^*$ in Theorem \ref{T24070901} is equal to 0, which means that $m(\rho)$ can be achieved for all $\rho>0$. Next, we aim to find the borderline to guarantee $\rho^*=0$. To this end, the cases that $\alpha_1=\alpha_2$ and $\alpha_1\neq \alpha_2$ have great differences. Actually the problem in the case $\alpha_1=\alpha_2$ is very similar to the single equation
\[
\begin{cases}
\Delta^2u+\alpha_{1}\Delta u+\lambda u=\beta |u|^{r-2} u & \text { in } \mathbb{R}^{N}, \\
\int_{\mathbb{R}^{N}} u^{2}\ud x=\rho^{2}.
\end{cases}
\]
We will in subsection \ref{2501121521} make use of the estimates of \cite{Fernandez22} to give a borderline for $\rho^*=0$. However, when $\alpha_1\neq \alpha_2$, the estimates are challenging and elusive.  In our paper, for $\alpha_1\neq  \alpha_2$, we obtain the condition   that can guarantee $\rho^*>0$, which is similar as the case $\alpha_1 =\alpha_2$, but it is difficult for us to conclude how  $\rho^*$ equals 0. Our result is the following.
\begin{theorem}\label{T2025011217}
Assume $\alpha_1,\,\alpha_2,\, \beta>0$ and $r\le 2+\frac{8}{N}$. Let $\rho^*$ be given in Theorem \ref{T24070901}. Then $\rho^*>0$  if
\begin{equation}\label{2411192127}
  \max\left\{4, 2+\frac{8}{N+1}\right\}\le r\le 2+\frac{8}{N}.
\end{equation}
Furthermore, when $\alpha_1=\alpha_2$, \eqref{2411192127} is also a necessary condition for $\rho^*>0$.
\end{theorem}
\begin{remark}
Precisely, under the assumptions that $r_1, r_2>1$ and $r_1+r_2\le 2+\frac{8}{N}$, Theorems \ref{T24070901} and \ref{T2025011217} say that when $\alpha_1=\alpha_2$, the nontrivial ground state exists for all $\rho>0$ if and only if
$$r_1+r_2<\min\left\{\max\left\{4, 2+\frac{8}{N+1}\right\}, 2+\frac{8}{N}\right\}.$$
\end{remark}

Now, we consider the case $r\in(\bar{r},2^{**})$. For the convenience of compactness, we will work in a subspace of $H^2(\mathbb{R}^N)\times H^2(\mathbb{R}^N)$, which consists of all radial functions, it is denoted by $H^{2}_{r}(\mathbb{R}^{N})\times H^{2}_{r}(\mathbb{R}^{N})$ where
\begin{equation}\label{2507102258}
\begin{aligned}
H^{2}_{r}(\mathbb{R}^{N}):=\left\{u\in H^{2}(\mathbb{R}^{N}):u \mbox{ is radially symmetric}\right\}.
\end{aligned}
\end{equation}
Observe that when the mass is small, the energy functional $I$ possesses a mountain pass structure. More precisely, when $\beta\rho^{r-2}<c^{*}(N,r)$, there exist $R_1>R_0>0$ such that for any $c\in (R_0, R_1)$, $I$ has a positive lower bound on $S_\rho$ with the restriction $\|\Delta u\|_2^2+\|\Delta v\|_2^2=c^2$, where
\begin{equation*}
\begin{aligned}
c^{*}(N,r):=
\frac{1}{2\mathcal{D}_1(r\gamma_{r}-1)}\left(\frac{r\gamma_{r}-2}
{(r\gamma_{r}-1)\max\{\alpha_{1},\alpha_{2}\}}\right)^{r\gamma_{r}-2},
\end{aligned}
\end{equation*}
$\gamma_{r}:=\frac{N(r-2)}{4r}$ and $\mathcal{D}_{1}$ is given in \eqref{2507061525}.
On the other hand, we can find two points $(u_0, v_0)$, $(u_1, v_1) \in S_\rho$ such that $I(u_0, v_0),\, I(u_1, v_1)<0$ with respectively $\|\Delta u_0\|_2^2+\|\Delta v_0\|_2^2<R_0^2$ and $\|\Delta u_1\|_2^2+\|\Delta v_1\|_2^2>R_1^2$. Thus, we can establish a mountain pass structure, precisely see Section \ref{2501190022}. Using Jeanjean's method in \cite{Jeanjean97}, one can obtain a Palais-Smale sequence $\{(u_n, v_n)\}$ approaching Poho\v{z}aev manifold, namely
\begin{equation}\label{2507102347}
2\int_{\mathbb{R}^{N}}(|\Delta u_n|^{2}+|\Delta v_n|^{2}) \ud x
-\alpha_{1}\int_{\mathbb{R}^{N}}|\nabla u_n|^{2} \ud x
-\alpha_{2}\int_{\mathbb{R}^{N}}|\nabla v_n|^{2} \ud x
-2\beta r\gamma_{r}\int_{\mathbb{R}^{N}}|u_n|^{r_{1}}|v_n|^{r_{2}} \ud x=o(1),
\end{equation}
and by virtue of Lagrange multiplier rule, we derive a sequence of Lagrange multipliers $\{\lambda_n\}$ corresponding to the Palais-Smale sequence $\{(u_n, v_n)\}$.  That is
\begin{align}\label{2507102333}
I'(u_n, v_n)+\lambda_n(u_n, v_n)\to 0\quad \mbox{ in } \left(H^2_r(\mathbb{R}^N)\times H^2_r(\mathbb{R}^N)\right)^*.
\end{align}

  As the usual process, it is not difficult to show that the sequence $(u_n, v_n)$ is bounded in $H_r^2(\mathbb{R}^N)\times H_r^2(\mathbb{R}^N)$, but whose compactness is very involved.  As previously emphasized, the main difficulty is that  we can not obtain the strong convergence $\nabla u_{n}\rightarrow\nabla u$ in $L^{2}(\mathbb{R}^{N})$ only from the weak convergence $u_{n} \rightharpoonup u$ in $H_r^{2}(\mathbb{R}^{N})$. Motivated by Luo and Yang \cite{Luo23}, we multiply \eqref{2507102333} by $(u_n, v_n)$, and combining \eqref{2507102347}, we can cancel the common terms $\int_{\mathbb{R}^{N}}|\nabla u_n|^{2} \ud x$ and $\int_{\mathbb{R}^{N}}|\nabla v_n|^{2} \ud x$. Thus the compactness could escape from dealing with the convergence of $\nabla u_n$ and $\nabla v_n$ in $L^2$. Moreover, if $(u_n, v_n)$ is not compact in $H^2_r(\mathbb{R}^N)\times H^2_r(\mathbb{R}^N)$, up to a subsequence, we can deduce  $\underset{n\to\infty}{\lim}\, \lambda_n\le \frac{\max\{\alpha^{2}_1,\alpha^{2}_2\}}{4}$, see \eqref{24033103}. Naturally, a key step to get the compactness is to obtain 
  \begin{equation}\label{2507110015}
 \underset{n\to\infty}{\lim}\, \lambda_n> \frac{\max\{\alpha^{2}_1,\alpha^{2}_2\}}{4}.
 \end{equation}
In Section \ref{2501190022}, we will prove that \eqref{2507110015} can holds if $\beta\rho^{r-2}<c_{*}(N,r)$ where
\begin{align*}
c_{*}(N,r)&:=\left\{
\begin{aligned}
&\frac{1}{2\mathcal{D}_1(r\gamma_{r}-1)}
\left(\frac{1-\gamma_{r}}{\gamma_{r}}\frac{4}{\max\{\alpha^{2}_1,\alpha^{2}_2\}}\right)^{\frac{r\gamma_{r}-2}{2}}\quad\quad\quad\text{ if }
\frac{1}{2}<\gamma_{r}<1,\\
&
\frac{1}{2\mathcal{D}_1(r\gamma_{r}-1)}
\left(\frac{r-2}{2(r\gamma_{r}-1)}\frac{4}{\max\{\alpha^{2}_1,\alpha^{2}_2\}}\right)^{\frac{r\gamma_{r}-2}{2}}~\quad\text{ if }
0<\gamma_{r}\leq\frac{1}{2}.
\end{aligned}
\right.
\end{align*}
Our result for $r\in (\bar{r}, 2^{**})$ can be stated as follows.
\begin{theorem}\label{T230828010}
Let $\alpha_{1}$, $\alpha_{2}$, $\beta>0$, $r_{1}$, $r_{2}>1$, $N\ge 2$ and $r\in(\bar{r},2^{**})$. Then for any $\rho>0$ satisfying
\begin{align*}
\beta \rho^{r-2}<\min\left\{c^{*}(N,r),c_{*}(N,r)\right\},
\end{align*}
problem \eqref{23082801}-\eqref{23101001} has a  nontrivial radial solution with mountain pass type for some $\lambda>\frac{\max\{\alpha^{2}_1,\alpha^{2}_2\}}{4}$.
\end{theorem}

\begin{remark}
Here we assume $N\ge 2$, because the compact embedding $H_r^2(\mathbb{R}^N)\subset L^p(\mathbb{R}^N)$ for $2<p<2^{**}$ only hold for $N\ge 2$.
\end{remark}

\medskip
This paper is organized as follows. In Section \ref{S2}, we present some useful lemmas. In Section \ref{S3}, we give the proof of Theorems \ref{T24070901} and \ref{T2025011217}. In Section \ref{2501190022}, we are dedicated to the proof of Theorem \ref{T230828010}. In the rest of this paper, unless otherwise specified, we always assume that $\alpha_{1}$, $\alpha_{2}$, $\beta$, $\rho>0$, $r_{1}$, $r_{2}>1$.
%

\section{Preliminaries}\label{S2}
\noindent

This section is devoted to some notations and preliminary results.

\begin{lemma}[Gagliardo-Nirenberg inequality \cite{Nirenberg59}]\label{L23083006}
 For $N \geq 1$ and $2<p<2^{**}$, there exists an optimal constant $C_{N,p}>0$ depending on $N$, $p$ such that
\begin{equation}\label{2411232346}
\|u\|_{p} \leq C_{N, p}\|u\|_{2}^{1-\gamma_{p}}\|\Delta u\|_{2}^{\gamma_{p}},\quad \forall u \in H^{2}(\mathbb{R}^{N}).
\end{equation}
Moreover, the extremal is attainable.
\end{lemma}

Applying classical Fourier transform and H\"older inequality, we can easily get the interpolation inequality
\begin{equation}\label{23092305}
\begin{aligned}
\|\nabla u\|^{2}_2\leq \|u\|_2\cdot\|\Delta u\|_{2},\quad \forall u \in H^{2}(\mathbb{R}^{N}).
\end{aligned}
\end{equation}
Using H\"{o}lder's inequality, Lemma \ref{L23083006} and Young's inequality, we have
\begin{align}\label{23091001}
\int|u|^{r_{1}}|v|^{r_{2}} \ud x&
\leq\left(\int\left|u\right|^{r}\right)^{\frac{r_{1}}{r}}
\left(\int\left|v\right|^{r}\right)^{\frac{r_{2}}{r}}\nonumber \\
& \leq C_{N,r}^r \left(\|u\|_2^{r_1}\|v\|_2^{r_2}\right)^{1-\gamma_{r}}
\left(\|\Delta u\|_2^{r_{1}}
\|\Delta v\|_2^{r_{2}}\right)^{\gamma_{r}}\\
& \leq \mathcal{D}_{1}\left(\|u\|_2^2+\|v\|_2^2\right)^{\frac{r\left(1-\gamma_{r}\right)}{2}}\left(\|\Delta u\|^{2}_{2}
+\|\Delta v\|^{2}_{2}\right)^{\frac{r \gamma_{r}}{2}}\nonumber,
\end{align}
where $\mathcal{D}_{1}=\left(\frac{r_1}{r}\right)^{r_1}\left(\frac{r_2}{r}\right)^{r_2} {C}_{N,r}^r$. In the last inequality above, we used the basic inequality
\begin{equation*}\label{2411241745}
a^{r_1}b^{r_2}\le \left(\frac{r_1}{r}\right)^{r_1}\left(\frac{r_2}{r}\right)^{r_2}\left(a^2+b^2\right)^\frac{r}{2}, \quad \forall a, b> 0,
\end{equation*}
where the equality holds if and only if $\frac{a^2}{b^2}=\frac{r_1}{r_2}$.
As a result, we obtain the following lemma.
\begin{lemma}\label{2412261437}
The equalities in \eqref{23091001} hold if and only if $v=\sqrt{\frac{r_2}{r_1}}u$ and $u$ is an extremal for {\it Gagliardo-Nirenberg} inequality given in \eqref{2411232346}.
\end{lemma}

By \eqref{23092305} and Cauchy inequality, we deduce that
\begin{equation}\label{24030301}
\begin{aligned}
\alpha_{1}\|\nabla u\|^{2}_2+\alpha_{2}\|\nabla v\|^{2}_2
\leq &\alpha_{1}\|u\|_2\|\Delta u\|_{2}+\alpha_{2}\|v\|_2\|\Delta v\|_{2}\\
\leq &\mathcal{D}_{2}\left(\|u\|_2^2+\|v\|_2^2\right)^{\frac{1}{2}}(\|\Delta u\|^{2}_{2}+\|\Delta v\|^{2}_{2})^{\frac{1}{2}},
\end{aligned}
\end{equation}
where $\mathcal{D}_{2}:=\max\{\alpha_{1},\alpha_{2}\}$.
\medskip

As the proof of \cite[Lemma 2.1]{Bonheure19} and \cite[Remark 3.10]{Fernandez22}, we give the Poho\v{z}aev identity for the problem \eqref{23082801} and omit it's proof.
\begin{lemma}\label{L23082901}
Assume that $(u,v)\in H^{2}(\mathbb{R}^{N})\times H^{2}(\mathbb{R}^{N})$ solves \eqref{23082801}, then the Poho\v{z}aev identity holds:
\begin{equation}\label{23083001}
\begin{aligned}
&\frac{N-4}{2} \left(\|\Delta u\|^{2}_2+\|\Delta v\|^{2}_2\right)
-\frac{N-2}{2}\left(\alpha_{1}\|\nabla u\|^{2}_2+\alpha_{2}\|\nabla v\|^{2}_2\right)
+\frac{N}{2} \lambda\left(\|u\|^{2}_2+\|v\|^{2}_2\right)\\
=&N\beta\int_{\mathbb{R}^{N}}|u|^{r_{1}}|v|^{r_{2}} \ud x.
\end{aligned}
\end{equation}
\end{lemma}
\begin{lemma}\label{L23101201}
Let $\alpha_1,\, \alpha_2,\, \beta>0$ and $r\in(2,2^{**}]$. When $\lambda\leq0$, then problem \eqref{23082801}-\eqref{23101001} admits no solutions.
\end{lemma}
\begin{proof}
Assume that $(u,v)$ is a solution of \eqref{23082801}-\eqref{23101001}. Multiplying the equation \eqref{23082801}  by $u$ and $v$ respectively, we get that
\begin{equation}\label{23062303}
\begin{aligned}
\|\Delta u\|^{2}_2+\|\Delta v\|^{2}_2
-\alpha_{1}\|\nabla u\|^{2}_2-\alpha_{2}\|\nabla v\|^{2}_2
+\lambda \|u\|^{2}_2+\lambda \|v\|^{2}_2
=\beta r\int_{\mathbb{R}^{N}}|u|^{r_{1}}|v|^{r_{2}} \ud x.
\end{aligned}
\end{equation}
By \eqref{23083001} and \eqref{23062303}, we have
\begin{equation}\label{23101202}
\begin{aligned}
-\alpha_{1}\|\nabla u\|^{2}_2-\alpha_{2}\|\nabla v\|^{2}_2
+2\lambda(\|u\|^{2}_2+\|v\|^{2}_2)
=\beta \left(N-r\frac{N-4}{2}\right)\int_{\mathbb{R}^{N}}|u|^{r_{1}}|v|^{r_{2}} \ud x.
\end{aligned}
\end{equation}
Since $\alpha_1,\, \alpha_2,\,\beta>0$ and $\lambda\leq0$, when $r\in(2,2^{**}]$ and by \eqref{23101202}, we can deduce that $\nabla u$, $\nabla v=0$ a.e. in $\mathbb{R}^{N}$, which contradicts $(u,v)\in S_\rho$.
\end{proof}

\section{Mass-subcritical and mass-critical case}\label{S3}
\noindent

In this section, we consider the existence of normalized solution for problem \eqref{23082801}-\eqref{23101001} with $r=r_{1}+r_{2}\in(2,\bar{r}]$. Initially, we show some important properties about $m(\rho)$ given in \eqref{2411071036}.
\begin{lemma}\label{l24070902}
Let $\rho>0$, $r\in(2,\bar{r}]$ and $\mathcal{D}_1$ be given in  \eqref{23091001}.

$\mathrm{(i)}$ $m(\rho)$ is finite if and only if either $r\in (2, \bar{r})$ or $r=\bar{r}$ and $\rho<\left(\frac{1}{2\mathcal{D}_1\beta}\right)^{\frac{N}{8}}$. Moreover, in case $r=\bar{r}$ and $\rho\ge\left(\frac{1}{2\mathcal{D}_1\beta}\right)^{\frac{N}{8}}$, we have $m(\rho)=-\infty$.

$\mathrm{(ii)}$ If either $r\in (2, \bar{r})$ or $r=\bar{r}$ and $\rho<\left(\frac{1}{2\mathcal{D}_1\beta}\right)^{\frac{N}{8}}$, then the map $\rho\mapsto m(\rho)$ is continuous.
\begin{proof}
(i) By \eqref{23091001} and \eqref{24030301}, we get that for $(u,v)\in S_{\rho}$,
\begin{equation}\label{23092806}
\begin{aligned}
I(u,v)=&\frac{1}{2}\|\Delta u\|^{2}_{2}+\frac{1}{2}\|\Delta v\|^{2}_{2}
-\frac{\alpha_{1}}{2}\|\nabla u\|^{2}_2
-\frac{\alpha_{2}}{2}\|\nabla v\|^{2}_2
-\beta\int_{\mathbb{R}^{N}}|u|^{r_{1}}|v|^{r_{2}} \ud x  \\
\geq& \frac{1}{2}\left(\|\Delta u\|^{2}_{2}+\|\Delta v\|^{2}_{2}\right)
-\frac{1}{2}\mathcal{D}_{2}\rho\left(\|\Delta u\|^{2}_{2}+\|\Delta v\|^{2}_{2}\right)^{\frac{1}{2}}\\
&-\mathcal{D}_{1}\beta\rho^{r\left(1-\gamma_{r}\right)}\left(\|\Delta u\|^{2}_{2}
+\|\Delta v\|^{2}_{2}\right)^{\frac{r \gamma_{r}}{2}}.
\end{aligned}
\end{equation}
When $r\in(2,\bar{r})$, we get that $r\gamma_{r}<2$ and $I$ is coercive on $S_\rho$,
which implies that $m(\rho)>-\infty$.

\smallskip

When $r=\bar{r}$ and $\mathcal{D}_{1}\beta \rho^{r\left(1-\gamma_{r}\right)}<\frac{1}{2}$, that is $\rho<\left(\frac{1}{2\mathcal{D}_1\beta}\right)^{\frac{N}{8}}$,  we have $r\gamma_{r}=2$ and $I$ is also coercive on $S_\rho$. Still $m(\rho)>-\infty$. If $\mathcal{D}_{1}\beta \rho^{r\left(1-\gamma_{r}\right)}\ge \frac{1}{2}$, we claim $I$ is unbounded from below on $S_\rho$. Let $U$ be an extremal of {\it Gagliardo-Nirenberg} inequality given in \eqref{2411232346} with $p=\bar{r}$ and $\|U\|_2^2=\frac{\rho^2 r_1}{r}$.
As we remarked in Lemma \ref{2412261437}, the equalities in \eqref{23091001} can hold if $u=U$ and $v=\sqrt{\frac{r_2}{r_1}}U$. Setting $U_t:=t^{\frac{N}{2}}U(t x)$, $u_t:=t^{\frac{N}{2}}u(t x)$ and $v_t:=t^{\frac{N}{2}}v(t x)$, then $\|U_t\|_2^2=\|U\|_2^2$ and
\[
\begin{aligned}
I(u_t, v_t)=&\left(\frac{1}{2}-\mathcal{D}_1\beta\rho^{r\left(1-\gamma_{r}\right)}\right)\frac{r}{r_1}\|\Delta U_t\|^{2}_{2}-\frac{\alpha_{1}}{2}\|\nabla U_t\|^{2}_{2}
-\frac{\alpha_{2}}{2}\frac{r_2}{r_1}\|\nabla U_t\|^{2}_{2}\\
\le &-\frac{\alpha_{1}t^2}{2}\|\nabla U\|^{2}_{2}
-\frac{\alpha_{2}t^2}{2}\frac{r_2}{r_1}\|\nabla U\|^{2}_{2}.
\end{aligned}
\]
Thus, letting $t\to\infty$ in the above, we get the desired conclusion.

\medskip
(ii) Let $\rho_{n}\rightarrow \rho$ as $n\rightarrow\infty$, and   in addition $\rho, \rho_n<\left(\frac{1}{2\mathcal{D}_1\beta}\right)^{\frac{N}{8}}$ if $r=\bar{r}$. There exists  $(u_{n},v_n)\in S_{\rho_n}$ such that
\begin{equation}\label{2411071126}
m(\rho_{n})\leq I(u_{n},v_n)<m\left(\rho_n\right)+\frac{1}{n}.
\end{equation}
By the coerciveness of $I$ on $S_\rho$, $\rho_n\to\rho$ and \eqref{2411071126}, we can easily deduce that $\left\{(u_{n},v_n)\right\}$ is bounded in $H^{2}(\mathbb{R}^{N})\times H^{2}(\mathbb{R}^{N})$. Let $\widetilde{u}_{n} :=\frac{\rho}{\rho_{n}} u_{n}$ and $\widetilde{v}_{n} :=\frac{\rho}{\rho_{n}} v_{n}$. Then $\{(\widetilde{u}_{n},\widetilde{v}_n)\}\subset S_\rho$ and
\begin{align*}
m(\rho)\leq& I(\widetilde{u}_{n},\widetilde{v}_n)\\
=&I(u_{n},v_n)
+\frac{1}{2}\left(\frac{\rho^{2}}{\rho_{n}^{2}}-1\right)
\int_{\mathbb{R}^{N}}(|\Delta u_{n}|^{2}+|\Delta v_{n}|^{2}) \ud x
-\frac{\alpha_{1}}{2}\left(\frac{\rho^{2}}{\rho_{n}^{2}}-1\right)\int_{\mathbb{R}^{N}}|\nabla u_{n}|^{2} \ud x\\
&-\frac{\alpha_{2}}{2}\left(\frac{\rho^{2}}{\rho_{n}^{2}}-1\right)\int_{\mathbb{R}^{N}}|\nabla v_{n}|^{2} \ud x
-\beta\left(\frac{\rho^{r}}{\rho_{n}^{r}}-1\right)\int_{\mathbb{R}^{N}}|u_{n}|^{r_{1}}|v_{n}|^{r_{2}} \ud x\\
=&I(u_{n},v_n)+o_{n}(1)
\end{align*}
due to $\rho_{n}\rightarrow \rho$. Thus
\begin{equation}\label{24070904}
m(\rho)\leq\liminf_{n\rightarrow\infty}m\left(\rho_{n}\right).
\end{equation}

On the other hand, define $\widetilde{w}^1_{n}:=\frac{\rho_{n}}{\rho}w^1_{n}$ and $\widetilde{w}^2_{n}:=\frac{\rho_{n}}{\rho}w^2_{n}$, where $\left\{(w^1_{n},w^2_{n})\right\}\subset S_\rho$ is a minimizing sequence for $m(\rho)$. By a similar
discussion as above,  $\left\{(w^1_{n},w^2_{n})\right\}$ is bounded in $H^{2}(\mathbb{R}^{N})\times H^{2}(\mathbb{R}^{N})$, and
\begin{equation*}
\begin{aligned}
m\left(\rho_{n}\right)\leq I\left(\widetilde{w}^1_{n},\widetilde{w}^2_{n}\right)=I\left(w^1_{n},w^2_{n}\right)+o_{n}(1)=m(\rho)+o_n(1).
\end{aligned}
\end{equation*}
Thus,
\begin{equation}\label{24070905}
\limsup_{n\rightarrow\infty}m\left(\rho_{n}\right)\leq m(\rho),
\end{equation}
Combining \eqref{24070904} with \eqref{24070905},
we obtain that $m(\rho)=\lim\limits_{n\rightarrow\infty}m(\rho_{n})$.
Thus (ii) holds.
\end{proof}
\end{lemma}

Let $m^{J}(\rho)$ be given in \eqref{2411111354}, which is well-defined for all $\rho>0$. he following result gives an accurate $m^{J}(\rho)$, and it is not attainable.
\begin{lemma}\label{24100101}
For any $\rho>0$, we have $m^{J}(\rho)=-\frac{\max\{\alpha_1^2, \alpha_2^2\}}{8}\rho^2$, which is never achieved.
\end{lemma}
\begin{proof}
By \eqref{2411111354}, we can rewrite
\begin{equation}\label{2411111427}
m^{J}(\rho)=\underset{\rho_1^2+\rho_2^2=\rho^2}{\inf}M(\rho_1, \rho_2),
\end{equation}
where
\[
M(\rho_1, \rho_2):=\underset{u, v\in H^2(\mathbb{R}^N),\,\|u\|_2^2=\rho_1^2,\, \|v\|_2^2=\rho_2^2}{\inf} J(u, v).
\]
It can easily see that
\[
\begin{aligned}
M(\rho_1, \rho_2)=&\underset{u\in H^2(\mathbb{R}^N),\,\|u\|_2^2=\rho_1^2}{\inf}\left(\frac{1}{2}\|\Delta u\|^{2}_{2}-\frac{\alpha_1}{2}\|\nabla u\|^{2}_2\right)\\
&+\underset{v\in H^2(\mathbb{R}^N),\,\|v\|_2^2=\rho_2^2}{\inf}\left(\frac{1}{2}\|\Delta v\|^{2}_{2}-\frac{\alpha_2}{2}\|\nabla v\|^{2}_2\right).
\end{aligned}
\]
As a result, it follows from \cite[Lemma 3.1]{Boussaid19} that
\begin{equation}\label{2411111428}
M(\rho_1, \rho_2)=-\frac{\alpha_1^2}{8}\rho_1^2-\frac{\alpha_2^2}{8}\rho_2^2.
\end{equation}
By \eqref{2411111427} and \eqref{2411111428}, we have
\[
m^{J}(\rho)=-\frac{\max\{\alpha_1^2, \alpha_2^2\}}{8}\rho^2.
\]

Now, we shall prove that $m^{J}(\rho)$ is never achieved.
When $\alpha_1\neq \alpha_2$,  we can simply suppose $\alpha_1>\alpha_2$. At this time, we have $m^{J}(\rho)=-\frac{\alpha_1^2}{8}\rho^2$. Suppose that there exists some $(u, v)\in  S_{\rho}$ such that $J(u, v)=-\frac{\alpha_1^2}{8}\rho^2$. We denote $\rho_1=\|u\|_2$ and $\rho_2=\|v\|_2 $, then
\[
-\frac{\alpha_1^2}{8}\rho^2=J(u, v)\ge M(\rho_1, \rho_2)=-\frac{\alpha_1^2}{8}\rho_1^2-\frac{\alpha_2^2}{8}\rho_2^2\ge -\frac{\alpha_1^2}{8}\rho^2,
\]
 hence $\rho_1=\rho$ and $\rho_2=0$. So
\[
\frac{1}{2}\|\Delta u\|^{2}_{2}-\frac{\alpha_1}{2}\|\nabla u\|^{2}_2=-\frac{\alpha_1^2}{8}\rho^2,
\]
and $u$ is a minimizer for
\[
M(\rho, 0)=\underset{u\in H^2(\mathbb{R}^N),\, \|u\|_2^2=\rho^2}{\inf} \left(\frac{1}{2}\|\Delta u\|^{2}_{2}-\frac{\alpha_1}{2}\|\nabla u\|^{2}_2\right).
\]
However we have already know that $M(\rho, 0)$ is never achieved, seeing \cite[Lemma 3.1]{Boussaid19} again. This contradiction tells us that $m^{J}(\rho)$ is also never achieved.

\smallskip
Suppose $\alpha_1=\alpha_2$ and $(u, v)\in  S_{\rho}$ is a minimizer of $m^{J}(\rho)$. With loss of generality, we assume $\rho_1=\|u\|_2\neq 0$. Thus, $u$ is a minimizer of $M(\rho_1, 0)$, which contradicts  \cite[Lemma 3.1]{Boussaid19}.
\end{proof}

With the above results, we can show the following  sub-additive argument.

\begin{lemma}\label{L24100202}
Assume either $r\in (2, \bar{r})$ or $r=\bar{r}$, $\rho<\left(\frac{1}{2\mathcal{D}_1\beta}\right)^{\frac{N}{8}}$. Suppose that  $m(\rho)$ is attained,  then

$\mathrm{(i)}$ for any $\theta>1$, we have $m(\theta\rho)<\theta^2 m(\rho)$;

$\mathrm{(ii)}$ for any $\rho_{1}>0$, we have $m\left(\left(\rho^2+\rho^2_{1}\right)^{\frac{1}{2}}\right)
< m(\rho)+m(\rho_{1})$, here $\rho_1<\left(\frac{1}{2\mathcal{D}_1\beta}\right)^{\frac{N}{8}}$ if $r=\bar{r}$.
\begin{proof}
Let  $(u,v)\in S_{\rho}$ be a global minimizer of $m(\rho)$. We first claim that $u, v\neq 0$. If otherwise, we can simply  assume $u\neq 0$ and $v=0$. Hence
\[
m^{J}(\rho)\le J(u, v)=I(u, v)=m(\rho).
\]
On the other hand, in view of the definition of $m^{J}(\rho)$ and $m(\rho)$, it is obvious that $m^{J}(\rho)\ge m(\rho)$. So $(u, v)$ is a minimizer of $m^{J}(\rho)$, which is a contradiction due to Lemma \ref{24100101}. Therefore the claim holds.

\smallskip
(i) For $\theta>1$, by $r\in(2,\bar{r}]$, we obtain
\begin{equation}\label{2411241444}
\begin{aligned}
I\left(\theta u,\theta v\right)
=\theta^{2}I\left(u,v\right)
+\beta\theta^{2}\left(1-\theta^{r-2}\right)
\int_{\mathbb{R}^{N}}|u|^{r_{1}}|v|^{r_{2}} \ud x
<\theta^{2}I\left(u,v\right)
=\theta^{2}m(\rho).
\end{aligned}
\end{equation}
Thus $m(\theta\rho)\leq I\left(\theta u,\theta v\right)<\theta^{2}I\left(u,v\right)
=\theta^{2}m(\rho)$. Hence (i) holds.

\smallskip
(ii) Without loss of generality, we assume that $\rho \geq \rho_{1}>0$. Then by (i),
\begin{equation}\label{2411241445}
\begin{aligned}
m\left(\left(\rho^2+\rho^2_{1}\right)^{1/2}\right)
=&m\left(\frac{\left(\rho^2+\rho^2_{1}\right)^{1/2}}{\rho}\rho\right)
<\frac{\rho^2+\rho^2_{1}}{\rho^2} m(\rho)\\
=&m(\rho)+\frac{\rho^{2}_{1}}{\rho^{2}} m(\rho)\leq m(\rho)+m(\rho_{1}),
\end{aligned}
\end{equation}
where we used $\frac{\rho^{2}_{1}}{\rho^{2}} m(\rho)\leq m(\rho_{1})$ in last inequality, whose proof is similar to (i), hence we omit it. Thus (ii) holds.
\end{proof}
\end{lemma}
\begin{lemma}\label{L24100205}
Assume that $\left\{(u_{n},v_n)\right\}$ is a bounded sequence in $H^{2}(\mathbb{R}^{N})\times H^{2}(\mathbb{R}^{N})$, and it satisfies $\lim\limits_{n \rightarrow \infty}\left(\|u_{n}\|^2_{2}+\|v_{n}\|^2_{2}\right)=\rho^2>0$. Let $\rho_{n}=\frac{\rho} {\left(\|u_{n}\|^2_{2}+\|v_{n}\|^2_{2}\right)^{1/2}}$, and $\widetilde{u}_{n}=\rho_{n} u_{n}$, $\widetilde{v}_{n}=\rho_{n} v_{n}$. Then the following holds:
\begin{equation*}
\begin{aligned}
(\widetilde{u}_{n},\widetilde{v}_{n})\in S_\rho, \quad \lim_{n \rightarrow \infty} \rho_{n}=1, \quad \lim_{n \rightarrow \infty}\left|I\left(\widetilde{u}_{n},\widetilde{v}_{n}\right)-I\left(u_{n},v_{n}\right)\right|=0.
\end{aligned}
\end{equation*}
\begin{proof}
Clearly, $(\widetilde{u}_{n},\widetilde{v}_{n})\in S_\rho$ and $\lim\limits_{n \rightarrow \infty}\rho_{n}=1$. So, using the boundedness of $\left\{(u_{n},v_n)\right\}$, we get that as $n \rightarrow \infty$,
\begin{align*}
I\left(\widetilde{u}_{n},\widetilde{v}_{n}\right)-I\left(u_{n},v_{n}\right)
=&\frac{\rho^2_n-1}{2}\int_{\mathbb{R}^{N}}(|\Delta u_{n}|^{2}+|\Delta v_{n}|^{2}) \ud x
-\frac{\alpha_{1}(\rho^2_n-1)}{2}\int_{\mathbb{R}^{N}}|\nabla u_{n}|^{2} \ud x\\
&-\frac{\alpha_{2}(\rho^2_n -1)}{2}\int_{\mathbb{R}^{N}}|\nabla v_{n}|^{2} \ud x
-\beta(\rho^r_n-1)\int_{\mathbb{R}^{N}}|u_{n}|^{r_{1}}|v_{n}|^{r_{2}} \ud x\\
\rightarrow &\,0.
\end{align*}
The proof is done.
\end{proof}
\end{lemma}

Next, we provide an alternative result for minimizing sequence. This will be used to conclude the existence of minimizers.

\begin{theorem}\label{T24100203}
Assume either $r\in (2, \bar{r})$ or $r=\bar{r}$, $\rho<\left(\frac{1}{2\mathcal{D}_1\beta}\right)^{\frac{N}{8}}$. Let $\{(u_n,v_n)\}\subset S_{\rho}$ be a minimizing sequence of $m(\rho)$ with $\rho>0$. Then one of the following holds:

$\mathrm{(i)}$ {\rm(}Vanishing{\rm )}
\begin{equation}\label{24100102}
\begin{aligned}
\limsup_{n\rightarrow\infty}\sup_{y\in\mathbb{R}^{N}}
\int_{B_{1}(y)}\left(|u_{n}|^{2}+|v_{n}|^{2}\right)\ud x=0.
\end{aligned}
\end{equation}

$\mathrm{(ii)}$ {\rm(}Compactness{\rm )} Up to a subsequence, there exist $(u,v)\in S_{\rho}$
and a sequence $\left\{y_{n}\right\} \subset \mathbb{R}^{N}$
such that
$$\left(u_{n}\left(\cdot-y_{n}\right),v_{n}\left(\cdot-y_{n}\right)\right) \rightarrow (u,v)\quad \mbox{in}\;\, H^{2}(\mathbb{R}^{N})\times H^{2}(\mathbb{R}^{N})$$
 as $n\rightarrow \infty$, and $(u,v)$ is a global minimizer.
\begin{proof}
Let $\{(u_n,v_n)\}\subset S_{\rho}$ be a minimizing sequence for $m(\rho)$, which does not satisfy (i). Thus,
\begin{equation}\label{2411080015}
\begin{aligned}
0<L:=\limsup_{n\rightarrow\infty}\sup_{y \in \mathbb{R}^{N}}
\int_{B_{1}(y)}\left(|u_{n}|^{2}+|v_{n}|^{2}\right) \ud x \leq \rho^{2},
\end{aligned}
\end{equation}
and there exists a sequence $\left\{y_{n}\right\} \subset \mathbb{R}^{N}$ such that up to a subsequence,
\begin{equation}\label{24100801}
\begin{aligned}
L=\lim_{n \rightarrow \infty}\int_{B_{1}(0)}
\left(|u_{n}(x-y_{n})|^{2}+|v_{n}(x-y_{n})|^{2}\right) \ud x.
\end{aligned}
\end{equation}
By the proof of Lemma \ref{l24070902}, we deduce that the sequence
$\left\{(u_{n},v_n)\right\}$ is bounded in $H^{2}(\mathbb{R}^{N})\times H^{2}(\mathbb{R}^{N})$.
So there exist $(u,v) \in H^{2}(\mathbb{R}^{N})\times H^{2}(\mathbb{R}^{N})$ and a renamed subsequence of $\{(u_n, v_n)\}$
such that
\begin{align*}
\left(u_{n}(\cdot-y_{n}),v_{n}(\cdot-y_{n})\right) \rightharpoonup& (u,v) \quad
\text { in } H^{2}(\mathbb{R}^{N})\times H^{2}(\mathbb{R}^{N}), \\
\left(u_{n}(\cdot-y_{n}),v_{n}(\cdot-y_{n})\right) \rightarrow& (u,v) \quad
\text { in } L^{p}_{loc}\left(\mathbb{R}^{N}\right)\times L^{p}_{loc}\left(\mathbb{R}^{N}\right)
 \text { for } 1\leq p<2^{**},\\
\left(u_{n}(\cdot-y_{n}),v_{n}(\cdot-y_{n})\right) \rightarrow& (u,v) \quad
\text { a.e.~in } \mathbb{R}^{N}\times \mathbb{R}^{N},
\end{align*}
and \eqref{2411080015}-\eqref{24100801} imply $(u,v)\neq (0,0)$. Now, let
\begin{equation*}
\begin{aligned}
\widetilde{u}_{n}:=u_{n}\left(\cdot-y_{n}\right)-u,\quad \widetilde{v}_{n}:=v_{n}\left(\cdot-y_{n}\right)-v.
\end{aligned}
\end{equation*}
By weak convergence and Brezis-Lieb type lemma \cite[Lemma 2.3]{CZ_TAMS2015}, we can obtain that
\begin{equation}\label{24100803}
\begin{aligned}
\left\|\Delta u_{n}\right\|_{2}^{2}=&\left\|\Delta\left(u+\widetilde{u}_{n}\right)\right\|_{2}^{2}
=\|\Delta u\|_{2}^{2}+\left\|\Delta \widetilde{u}_{n}\right\|_{2}^{2}+o_{n}(1), \\
\left\|\Delta v_{n}\right\|_{2}^{2}=&\left\|\Delta\left(v+\widetilde{v}_{n}\right)\right\|_{2}^{2}
=\|\Delta v\|_{2}^{2}+\left\|\Delta \widetilde{v}_{n}\right\|_{2}^{2}+o_{n}(1), \\
\left\|\nabla u_{n}\right\|_{2}^{2}=&\left\|\nabla\left(u+\widetilde{u}_{n}\right)\right\|_{2}^{2}
=\|\nabla u\|_{2}^{2}+\left\|\nabla \widetilde{u}_{n}\right\|_{2}^{2}+o_{n}(1), \\
\left\|\nabla v_{n}\right\|_{2}^{2}=&\left\|\nabla\left(v+\widetilde{v}_{n}\right)\right\|_{2}^{2}
=\|\nabla v\|_{2}^{2}+\left\|\nabla \widetilde{v}_{n}\right\|_{2}^{2}+o_{n}(1),
\end{aligned}
\end{equation}
\begin{equation}\label{24100804}
\begin{aligned}
\int_{\mathbb{R}^{N}}\left(|{u}_{n}|^{2}+|{v}_{n}|^{2}\right)\ud x
=&\int_{\mathbb{R}^{N}}\left(|u+\widetilde{u}_{n}|^{2}+|v+\widetilde{v}_{n}|^{2}\right)\ud x\\
=&\int_{\mathbb{R}^{N}}\left(|{u}|^{2}+|{v}|^{2}\right)\ud x
+\int_{\mathbb{R}^{N}}\left(|\widetilde{u}_{n}|^{2}+|\widetilde{v}_{n}|^{2}\right)\ud x+o_{n}(1)
\end{aligned}
\end{equation}
and
\begin{equation*}
\begin{aligned}
\int_{\mathbb{R}^{N}}|u_{n}|^{r_{1}}|v_{n}|^{r_{2}} \ud x
=&\int_{\mathbb{R}^{N}}|u+\widetilde{u}_{n}|^{r_{1}}|v+\widetilde{v}_{n}|^{r_{2}} \ud x\\
=&\int_{\mathbb{R}^{N}}|u|^{r_{1}}|v|^{r_{2}} \ud x
+\int_{\mathbb{R}^{N}}|\widetilde{u}_{n}|^{r_{1}}|\widetilde{v}_{n}|^{r_{2}} \ud x
+o_{n}(1).
\end{aligned}
\end{equation*}
Hence
\begin{equation}\label{24100805}
\begin{aligned}
I(u_n,v_n)=I(u_n\left(x-y_{n}\right),v_n\left(x-y_{n}\right))
=I(u,v)+I(\widetilde{u}_{n},\widetilde{v}_{n})+o_{n}(1).
\end{aligned}
\end{equation}

{\it Claim.} $\int_{\mathbb{R}^{N}}\left(|\widetilde{u}_{n}|^{2}+|\widetilde{v}_{n}|^{2}\right)\ud x\rightarrow 0$ as $n \rightarrow \infty$. That is, $\int_{\mathbb{R}^{N}}\left(|{u}|^{2}+|{v}|^{2}\right)\ud x=\rho^2$.

\smallskip
In fact, let $\rho^2_{1}=\int_{\mathbb{R}^{N}}\left(|{u}|^{2}+|{v}|^{2}\right)\ud x>0$.
By \eqref{24100804}, if we get $\rho_{1}=\rho$, the claim follows.
Assume $\rho_{1}<\rho$, and define
\begin{equation*}
\begin{aligned}
\hat{u}_{n}=\left(\frac{\rho^2-\rho^2_{1}}
{\left\|\widetilde{u}_{n}\right\|^2_{2}+\left\|\widetilde{v}_{n}\right\|^2_{2}}\right)^{1/2} \widetilde{u}_{n},
\quad \hat{v}_{n}=\left(\frac{\rho^2-\rho^2_{1}}
{\left\|\widetilde{u}_{n}\right\|^2_{2}+\left\|\widetilde{v}_{n}\right\|^2_{2}}\right)^{1/2} \widetilde{v}_{n}.
\end{aligned}
\end{equation*}
By \eqref{24100805} and Lemma \ref{L24100205}, it follows that
\begin{equation*}
\begin{aligned}
I(u_n,v_n)=&I(u,v)+I(\widetilde{u}_{n},\widetilde{v}_{n})+o_{n}(1)\\
=&I(u,v)+I(\hat{u}_{n},\hat{v}_{n})+o_{n}(1)\\
\geq& I(u,v)+m\left((\rho^2-\rho^2_{1})^{1/2}\right)+o_{n}(1).
\end{aligned}
\end{equation*}
Hence, similar to \eqref{2411241444} and \eqref{2411241445}, we have
\begin{equation}\label{24100806}
\begin{aligned}
m(\rho) \geq I(u,v)+m\left((\rho^2-\rho^2_{1})^{{1}/{2}}\right)
\geq m\left(\rho_{1}\right)+m\left((\rho^2-\rho^2_{1})^{1/2}\right) \geq m(\rho).
\end{aligned}
\end{equation}
Thus $I(u,v)=m\left(\rho_{1}\right)$. That is, $(u,v)$ is global minimizer with respect to $\rho_{1}$.
Using Lemma \ref{L24100202}-(ii), we deduce the strict inequality
\begin{equation*}
\begin{aligned}
m(\rho)< m\left(\rho_{1}\right)+m\left((\rho^2-\rho^2_{1})^{1/2}\right),
\end{aligned}
\end{equation*}
which contradicts \eqref{24100806}. So $\rho_{1}=\rho$, hence we complete the proof of the claim.

\smallskip
Since $\{\left(\widetilde{u}_{n},\widetilde{v}_{n}\right)\}$ is a bounded sequence in $H^{2}(\mathbb{R}^{N})\times H^{2}(\mathbb{R}^{N})$ and using the above claim, it follows from \eqref{23092305} and \eqref{23091001} respectively that $\left\|\nabla \widetilde{u}_{n}\right\|_{2}^{2} \rightarrow 0$, $\left\|\nabla \widetilde{v}_{n}\right\|_{2}^{2} \rightarrow 0$ and $\int_{\mathbb{R}^{N}}|\widetilde{u}_{n}|^{r_{1}}|\widetilde{v}_{n}|^{r_{2}} \ud x \rightarrow 0$ as $n \rightarrow \infty$. Thus,
\begin{equation}\label{24100807}
\begin{aligned}
\liminf_{n \rightarrow \infty} I(\widetilde{u}_{n},\widetilde{v}_{n})
=\liminf_{n \rightarrow \infty} \frac{1}{2}
\left(\|\Delta \widetilde{u}_{n}\|_{2}^{2}+\|\Delta \widetilde{v}_{n}\|_{2}^{2}\right) \geq 0.
\end{aligned}
\end{equation}

On the other hand, since $\int_{\mathbb{R}^{N}}\left(|{u}|^{2}+|{v}|^{2}\right)\ud x=\rho^2$,
we deduce from \eqref{24100805} that
\begin{equation*}
\begin{aligned}
I(u_n,v_n)=I(u,v)+I(\widetilde{u}_{n},\widetilde{v}_{n})+o_{n}(1)
 \geq m(\rho)+I(\widetilde{u}_{n},\widetilde{v}_{n})+o_{n}(1),
\end{aligned}
\end{equation*}
and so
\begin{equation}\label{24100808}
\begin{aligned}
\limsup_{n \rightarrow \infty}I(\widetilde{u}_{n},\widetilde{v}_{n}) \leq 0.
\end{aligned}
\end{equation}
From \eqref{24100807} and \eqref{24100808}, we get that
\begin{equation*}
\begin{aligned}
\left\|\Delta \widetilde{u}_{n}\right\|_{2}^{2} \rightarrow 0,
\quad\left\|\Delta \widetilde{v}_{n}\right\|_{2}^{2} \rightarrow 0.
\end{aligned}
\end{equation*}
So, by \eqref{24100803}, $\left(u_{n}\left(\cdot-y_{n}\right),v_{n}\left(\cdot-y_{n}\right)\right) \rightarrow (u,v)$ in $H^{2}(\mathbb{R}^{N})\times H^{2}(\mathbb{R}^{N})$ as $n\rightarrow \infty$. Hence $I(u,v)=m(\rho)$, and by Lemma \ref{L24100202}, $(u,v)\in S_\rho$ is a global minimizer.
\end{proof}
\end{theorem}

With the above alternative argument, we have the following criterion to derive the existence of minimizer for $m(\rho)$.

\begin{proposition}\label{2411122357}
Assume either $r\in (2, \bar{r})$ or $r=\bar{r}$, $\rho<\left(\frac{1}{2\mathcal{D}_1\beta}\right)^{\frac{N}{8}}$. For $\rho>0$, let $m(\rho)$ and $m^{J}(\rho)$ be given respectively in \eqref{2411071036} and \eqref{2411111354}. Then $m(\rho)$ is achieved if $m(\rho)<m^{J}(\rho)$.
\end{proposition}
\begin{proof}
Suppose that $m(\rho)<m^{J}(\rho)$. Let $\{(u_n, v_n)\}$ be a minimizing sequence of $m(\rho)$. If the sequence $\{(u_n, v_n)\}$ vanishes, according to Lions' lemma \cite[Lemma I.1]{Lions842}, we have
\begin{align*}
\left(u_{n}, v_{n}\right) \rightharpoonup& (0, 0) \
\text { in } H^{2}(\mathbb{R}^{N})\times H^{2}(\mathbb{R}^{N}), \\
\left(u_{n},v_{n}\right) \rightarrow& (0, 0) \
\text { in } L^{p}\left(\mathbb{R}^{N}\right)\times L^{p}(\mathbb{R}^{N}),
\quad \text { for } 2< p<2^{**}.
\end{align*}
Thus
\[
\begin{aligned}
m(\rho)&=I(u_n, v_n)+o_n(1)\\
&=\frac{1}{2}\left(\|\Delta u_n\|^{2}_{2}+\|\Delta v_n\|^{2}_{2}\right)
-\frac{1}{2}\alpha_{1}\|\nabla u_n\|^{2}_2-\frac12 \alpha_{2}
\|\nabla v_n\|^{2}_2+o_n(1)\\
&= J(u_n, v_n)+o_n(1)\\
&\ge m^{J}(\rho)+o_n(1).
\end{aligned}
\]
This is a contradiction with $m(\rho)<m^{J}(\rho)$, so $\{(u_n, v_n)\}$  does not vanish. Using Theorem \ref{T24100203}, $\{(u_n, v_n)\}$ has a convergent subsequence and $m(\rho)$ is achieved.
\end{proof}

\subsection{Comparison between $m(\rho)$ and $m^J(\rho)$}\label{2501182319}
In this subsection, we will explore what conditions could guarantee $m(\rho)<m^J(\rho)$. Now we assume $\alpha_1\ge \alpha_2$, hence $m^J(\rho)=-\frac{\alpha_1^2}{8}\rho^2$. For any $(u, v)\in S_\rho$, we have
\begin{equation}\label{2411130006}
\begin{aligned}
H(u, v):=&I(u, v)-m^J(\rho)\\
=&\frac12 \left\|\left(\Delta + \frac{\alpha_1}{2}\right)u\right\|_2^2+\frac12 \left\|\left(\Delta + \frac{\alpha_2}{2}\right)v\right\|_2^2+\frac{\alpha_1^2-\alpha_2^2}{8}\|v\|_2^2-\beta\int_{\mathbb{R}^{N}}|u|^{r_{1}}|v|^{r_{2}} \ud x\\
=&\frac12 \left(\left\|\left(\Delta + \frac{\alpha_1}{2}\right)u\right\|_2^2+ \left\|\left(\Delta + \frac{\alpha_2}{2}\right)v\right\|_2^2+\frac{\alpha_1^2-\alpha_2^2}{4}\|v\|_2^2\right)\Big(1-2\beta\rho^{r-2}Q(u, v)\Big),
\end{aligned}
\end{equation}
where $Q$ is defined in \eqref{2411192136}. Let us denote
\begin{equation}\label{2411192130}
R:=\sup\, \left\{Q(u, v): (u, v)\in (H^2(\mathbb{R}^{N})\times H^2(\mathbb{R}^{N}))\backslash \{(0, 0)\}\right\}.
\end{equation}
Next, if $R<\infty$, we denote
\begin{equation}\label{2411241620}
\rho^*:=\left(\frac{1}{2\beta R}\right)^{\frac{1}{r-2}}.
\end{equation}
We will see in Proposition \ref{2412261421} below that $\rho^*$ is the dichotomy parameter for whether $m(\rho)<m^J(\rho)$. Note that in the mass-critical case $r=\bar{r}$, Lemma \ref{l24070902} tells us that $m(\rho)>-\infty$ only holds for $\rho<\left(\frac{1}{2\mathcal{D}_1\beta}\right)^{\frac{N}{8}}$. Hence it is important for us whether $\rho^*<\left(\frac{1}{2\mathcal{D}_1\beta}\right)^{\frac{N}{8}}$.
\begin{lemma}\label{2412261533}
If $r=\bar{r}$, there holds that $\rho^*<\left(\frac{1}{2\mathcal{D}_1\beta}\right)^{\frac{N}{8}}$.
\end{lemma}
\begin{proof}
In view of \eqref{2411241620}, the inequality $\rho^*<\left(\frac{1}{2\mathcal{D}_1\beta}\right)^{\frac{N}{8}}$ is equivalent to $R>\mathcal{D}_1$.
To this end, it suffices to show that there exists some $(u_0, v_0)$ such that
\[
Q(u_0, v_0)>\mathcal{D}_1.
\]
For that, let $U$ be an extremal of {\it Gagliardo-Nirenberg} inequality given in \eqref{2411232346} with $p=\bar{r}$. Now we take  $U_t:=U(t x)$,  and $u_t:=U_t$ and $v_t:=\sqrt{\frac{r_2}{r_1}}U_t$. Thus
\[
Q(u_t, v_t)=\frac{{(\frac{r_2}{r_1})}^{\frac{r_2}{2}}\int_{\mathbb{R}^{N}}|U_t|^{r} \ud x}{\left(1+\frac{r_2}{r_1}\right)^{\frac{r}{2}-1}\left(\left\|\left(\Delta + \frac{\alpha_1}{2}\right)U_t\right\|_2^2+ \frac{r_2}{r_1}\left\|\left(\Delta + \frac{\alpha_2}{2}\right)U_t\right\|_2^2+\frac{\alpha_1^2-\alpha_2^2}{4}\frac{r_2}{r_1}\|U_t\|_2^2\right)\|U_t\|_2^{r-2}}.
\]
By \eqref{23091001} and Lemma \ref{2412261437}, we have
\[
\left(\frac{r_2}{r_1}\right)^{\frac{r_2}{2}}\|U_t\|_r^r=\mathcal{D}_{1}\left(1+\frac{r_2}{r_1}\right)^{\frac{r}{2}}\|U_t\|_2^{r-r\gamma_r}\|\Delta U_t\|_2^{r\gamma_r}.
\]
Consequently, we obtain
\[
\begin{aligned}
Q(u_t, v_t)=&\frac{\mathcal{D}_{1}\left(1+\frac{r_2}{r_1}\right)\|\Delta U_t\|_2^{r\gamma_r}}{\left\|\left(\Delta + \frac{\alpha_1}{2}\right)U_t\right\|_2^2+ \frac{r_2}{r_1}\left\|\left(\Delta + \frac{\alpha_2}{2}\right)U_t\right\|_2^2+\frac{\alpha_1^2-\alpha_2^2}{4}\frac{r_2}{r_1}\|U_t\|_2^2}\\
=&\frac{\mathcal{D}_{1}\left(1+\frac{r_2}{r_1}\right)t^{4-N}\|\Delta U\|_2^{2}}{\left(1+\frac{r_2}{r_1}\right)t^{4-N}\|\Delta U\|_2^{2}-\left(\alpha_1+\frac{r_2}{r_1}\alpha_2\right)t^{2-N}\|\nabla U\|_2^2+ \left(\frac{\alpha_1^2}{4}+\frac{r_2\alpha_2^2}{4r_1}+\frac{\alpha_1^2-\alpha_2^2}{4}\right)t^{-N}\|U\|_2^2}.\\
\end{aligned}
\]
It is easy to get $Q(u_t, v_t)>\mathcal{D}_1$ as $t$ is large enough. The proof is done.
\end{proof}

\begin{proposition}\label{2412261421}
Assume either $r\in (2, \bar{r})$ or $r=\bar{r}$, $\rho<\left(\frac{1}{2\mathcal{D}_1\beta}\right)^{\frac{N}{8}}$. Then we have
\begin{itemize}
\item[\rm (i)] when $R<\infty$,  there hold that $m(\rho)=m^{J}(\rho)$ if $\rho\in (0,\rho^*]$, and $m(\rho)<m^{J}(\rho)$ if $\rho>\rho^*$. Moreover, when $\rho\in (0, \rho^*)$, $m(\rho)$ is never achieved;
\item[\rm (ii)]when $R=\infty$, there hold that  $m(\rho)<m^J(\rho)$ for all $\rho>0$.
\end{itemize}
\end{proposition}
\begin{proof}
{\rm (i)} If $\rho\in (0,\rho^*]$, for any $(u, v)\in S_{\rho}$, we obtain
\[
0\le 1-2\beta\rho^{r-2}R\le 1-2\beta\rho^{r-2}Q(u, v).
\]
Hence it follows from \eqref{2411130006} that
\begin{equation}
\underset{(u, v)\in S_{\rho}}{\inf}H(u, v)\ge 0,
\end{equation}
so $m(\rho)=m^{J}(\rho)$.

Next, we claim that $m(\rho)$ is never achieved for all $\rho\in (0, \rho^*)$. We have already known that $m(\rho)=-\frac{\alpha_1^2}{8}\rho^2$ for all $\rho\in (0, \rho^*)$, so $m$ is differentiable and
\begin{equation}\label{2411201553}
m'(\rho)=-\frac{\alpha_1^2}{4}\rho.
\end{equation}
Now we assume by contradiction that $(u, v)$ is a minimizer for $m(\rho)$. Observe that $u\neq 0$ and $v\neq 0$. If otherwise, $(u, v)$ is also a minimizer for $m^J(\rho)$, that is a contradiction by Lemma \ref{24100101}.

Taking $t>0$ such that $(1+t)\rho<\rho^*$, we have
\[
I((1+t)u, (1+t)v)\ge m((1+t)\rho),
\]
which together with $I(u, v)=m(\rho)$ implies that
\begin{equation}\label{2411201546}
\begin{aligned}
m'(\rho)\rho=& \underset{t\to 0^+}{\lim}\frac{m((1+t)\rho)-m(\rho)}{t}\\
\le& \underset{t\to 0^+}{\lim\inf} \frac{I((1+t)u, (1+t)v)-I(u, v)}{t}\\
=&\langle I'(u, v), (u, v)\rangle.
\end{aligned}
\end{equation}
On the other hand, for $0<t<1$, we have
\[
I((1-t)u, (1-t)v)\ge m((1-t)\rho),
\]
then
\begin{equation}\label{2411201547}
\begin{aligned}
m'(\rho)\rho=& \underset{t\to 0^+}{\lim}\frac{m((1-t)\rho)-m(\rho)}{-t}\\
\ge& \underset{t\to 0^+}{\lim\sup} \frac{I((1-t)u, (1-t)v)-I(u, v)}{-t}\\
=&\langle I'(u, v), (u, v)\rangle.
\end{aligned}
\end{equation}
From \eqref{2411201553}-\eqref{2411201547}, it follows that
\begin{equation}\label{2411201608}
\langle I'(u, v), (u, v)\rangle=-\frac{\alpha_1^2}{4}\rho^2.
\end{equation}
According to the Lagrange multiplier principle, there exists a $\lambda_\rho\in\mathbb{R}$ such that
\[
I'(u, v)=\lambda_\rho(u, v) \quad \mbox{ in}\; \left(H^{2}(\mathbb{R}^N)\times H^{2}(\mathbb{R}^N)\right)^{-1}.
\]
This combined with \eqref{2411201608} gives
\[
\lambda_\rho=-\frac{\alpha_1^2}{4}.
\]
Since $(u, v)$ is a solution to \eqref{23082801} with $\lambda=-\lambda_\rho$, we get
\[
\|\Delta u\|^{2}_2+\|\Delta v\|^{2}_2
-\alpha_{1}\|\nabla u\|^{2}_2-\alpha_{2}\|\nabla v\|^{2}_2
+\frac{\alpha_1^2\rho^2}{4}
=\beta r\int_{\mathbb{R}^{N}}|u|^{r_{1}}|v|^{r_{2}} \ud x.
\]
Thus
\[
-\frac{\alpha_1^2\rho^2}{8}=m(\rho)=I(u, v)=-\frac{\alpha_1^2\rho^2}{8}+\left(\frac{r}{2}-1\right)\beta\int_{\mathbb{R}^{N}}|u|^{r_{1}}|v|^{r_{2}} \ud x.
\]
This is a contradiction since $\frac{r}{2}>1$ and $u\neq 0$, $v\neq 0$.

\medskip
{\rm (ii)} If $R=\infty$, there must be some $(u, v)\in S_\rho$ such that the last term of \eqref{2411130006} is negative, which means that
\begin{equation}\label{2411130054}
\underset{(u, v)\in S_{\rho}}{\inf}H(u, v)<0,
\end{equation}
that is $m(\rho)<m^J(\rho)$.
\end{proof}
\begin{proof}[\bf Proof~of~Theorem~\ref{T24070901} completed] Let $R$ be given as in \eqref{2411192130}, and
\[
\begin{aligned}
\rho^*=
\begin{cases}
\left(\frac{1}{2\beta R}\right)^{\frac{1}{r-2}} \quad & \mbox{if}\; R<\infty;\\
0 \quad &\mbox{if}\; R=\infty.
\end{cases}
\end{aligned}
\]
By Lemma \ref{2412261533} and Proposition \ref{2412261421}, we have $m(\rho)<m^J(\rho)$ if
\[
\begin{aligned}
\rho\in
\begin{cases}
 (\rho^*, \infty)\quad &\mbox{if}\;\; r<\bar{r};\\
\left(\rho^*, \left(\frac{1}{2\mathcal{D}_1\beta}\right)^{\frac{N}{8}}\right)\quad &\mbox{if}\;\; r=\bar{r}.
\end{cases}
\end{aligned}
\]
Using Proposition \ref{2411122357}, $m(\rho)$ is achieved. By Proposition \ref{2412261421} again, when $\rho<\rho^*$, $m(\rho)$ is never achieved. Therefore we can complete the proof.
\end{proof}
\subsection{Estimates on $R$}\label{2501121521}
According to Theorem \ref{T24070901},  $m(\rho)$ is achieved for all $\rho>0$ if $\rho^*=0$. Thus, it is very important for us to know what conditions can guarantee $\rho^*=0$ or $\rho^*>0$. For that, we will discuss the value of $R$ by two cases: $\alpha_1=\alpha_2$ and $\alpha_1\neq\alpha_2$.
\subsection*{(i) Case $\alpha_1=\alpha_2$}

\begin{proposition}\label{2501171415}
Assume that $\alpha_1=\alpha_2$. Let $R$ be given in \eqref{2411192130}, then $R$ is finite if and only if  \eqref{2411192127} holds.
\end{proposition}
\begin{proof}
For simplicity, we set $\alpha:=\alpha_1=\alpha_2$. By H\"older inequality and Young's inequality, we can calculate that
\begin{equation}\label{2412261602}
\begin{aligned}
Q(u, v)\le& \frac{\|u\|_r^{r_1}\|v\|_r^{r_2}}{\left(\left\|\left(\Delta + \frac{\alpha}{2}\right)u\right\|_2^2+ \left\|\left(\Delta + \frac{\alpha}{2}\right)v\right\|_2^2\right)\left(\|u\|_2^2+\|v\|_2^2\right)^{\frac{r}{2}-1}}\\
\le& \frac{\frac{r_1}{r}\|u\|_r^{r}+\frac{r_2}{r}\|v\|_r^{r}}{\left(\left\|\left(\Delta + \frac{\alpha}{2}\right)u\right\|_2^2+ \left\|\left(\Delta + \frac{\alpha}{2}\right)v\right\|_2^2\right)\left(\|u\|_2^2+\|v\|_2^2\right)^{\frac{r}{2}-1}}\\
\le&\frac{\frac{r_1}{r}\|u\|_r^{r}+\frac{r_2}{r}\|v\|_r^{r}}{\left(\left\|\left(\Delta + \frac{\alpha}{2}\right)u\right\|_2^2+ \left\|\left(\Delta + \frac{\alpha}{2}\right)v\right\|_2^2\right)\cdot \frac12\left(\|u\|_2^{r-2}+\|v\|_2^{r-2}\right)}\\
\le&\frac{\frac{r_1}{r}\|u\|_r^{r}+\frac{r_2}{r}\|v\|_r^{r}}{\frac12\left(\left\|\left(\Delta + \frac{\alpha}{2}\right)u\right\|_2^2\|u\|_2^{r-2}+ \left\|\left(\Delta + \frac{\alpha}{2}\right)v\right\|_2^2\|v\|_2^{r-2}\right)}\\
\le &\frac{2\max\{r_1, r_2\}}{r}\sup_{u\in H^2(\mathbb{R}^{N})\backslash \{0\}}\frac{\|u\|_r^r}{\left\|\left(\Delta + \frac{\alpha}{2}\right)u\right\|_2^2\|u\|_2^{r-2}}.
\end{aligned}
\end{equation}
Note that in the last inequality above, we used the basic inequality
\begin{equation}\label{2411160031}
\frac{a+b}{c+d}\le \max\left\{\frac{a}{c}, \frac{b}{d}\right\},\quad \forall a, b, c, d>0.
\end{equation}
Therefore, provided that
\begin{equation*}
\widetilde{R}:=\sup_{u\in H^2(\mathbb{R}^{N})\backslash \{0\}}\frac{\|u\|_r^r}{\left\|\left(\Delta + \frac{\alpha}{2}\right)u\right\|_2^2\|u\|_2^{r-2}}<\infty,
\end{equation*}
then we can get $R<\infty$.

By the scaling transformation $w=u\left(\sqrt{\frac{\alpha}{2}}x\right)$, we have
\begin{equation*}
\begin{aligned}
\|\Delta w\|_2^2=\left({\frac{2}{\alpha}}\right)^{\frac{N-4}{2}}\|\Delta u\|_2^2,\quad \|\nabla w\|_2^2=\left({\frac{2}{\alpha}}\right)^{\frac{N-2}{2}}\|\nabla u\|_2^2,\quad\| w\|_r^r=\left({\frac{2}{\alpha}}\right)^{\frac{N}{2}}\| u\|_r^r.
\end{aligned}
\end{equation*}
Consequently,
\begin{equation}\label{2411192049}
\sup_{w\in H^2(\mathbb{R}^{N})\backslash \{0\}}\frac{\|w\|_r^r}{\left\|\left(\Delta + \frac{\alpha}{2}\right)w\right\|_2^2\|w\|_2^{r-2}}=\left(\frac{2}{\alpha}\right)^{2-\frac{(r-2)N}{4}}\sup_{u\in H^2(\mathbb{R}^{N})\backslash \{0\}}\frac{\|u\|_r^r}{\left\|\left(\Delta + 1\right)u\right\|_2^2\|u\|_2^{r-2}}.
\end{equation}
By virtue of \cite[Theorem 1.1]{Fernandez22}, we have that \eqref{2411192049} is finite if and only if \eqref{2411192127} holds. Combined with \eqref{2412261602}, therefore \eqref{2411192127}  can also guarantee $R<\infty$.

In the following, we shall prove that $R=\infty$ if \eqref{2411192127} does not hold. In fact, we can take $u=v$ in \eqref{2411192136}, that is
\[
Q(u, u)= \frac{\int_{\mathbb{R}^{N}}|u|^{r} \ud x}{2^{\frac{r}{2}}\left\|\left(\Delta + \frac{\alpha}{2}\right)u\right\|_2^2\|u\|_2^{r-2}},
\]
which comes back to the analysis of \eqref{2411192049}. The proof is done.
\end{proof}
\subsection*{(ii) Case $\alpha_1\neq \alpha_2$}
Without confusions, we assume $\alpha_2<\alpha_1$. Similar to the case $\alpha_1=\alpha_2$, here we still compute the supremum of $Q(u, v)$. By H\"older inequality, Young's inequality and \eqref{2411160031}, we can get
\[
\begin{aligned}
Q(u, v)\le& \frac{\|u\|_r^{r_1}\|v\|_r^{r_2}}{\left(\left\|\left(\Delta + \frac{\alpha_1}{2}\right)u\right\|_2^2+ \left\|\left(\Delta + \frac{\alpha_2}{2}\right)v\right\|_2^2\right)\left(\|u\|_2^2+\|v\|_2^2\right)^{\frac{r}{2}-1}}\\
\le&\frac{\frac{r_1}{r}\|u\|_r^{r}+\frac{r_2}{r}\|v\|_r^{r}}{\left(\left\|\left(\Delta + \frac{\alpha_1}{2}\right)u\right\|_2^2+ \left\|\left(\Delta + \frac{\alpha_2}{2}\right)v\right\|_2^2\right)\cdot \frac12\left(\|u\|_2^{r-2}+\|v\|_2^{r-2}\right)}\\
\le&\frac{\frac{r_1}{r}\|u\|_r^{r}+\frac{r_2}{r}\|v\|_r^{r}}{\frac12\left(\left\|\left(\Delta + \frac{\alpha_1}{2}\right)u\right\|_2^2\|u\|_2^{r-2}+ \left\|\left(\Delta + \frac{\alpha_2}{2}\right)v\right\|_2^2\|v\|_2^{r-2}\right)}\\
\le &\frac{2\max\{r_1, r_2\}}{r}\max_{j=1, 2}\left\{\sup_{u\in H^2(\mathbb{R}^{N})\backslash \{0\}}\frac{\|u\|_r^r}{\left\|\left(\Delta + \frac{\alpha_j}{2}\right)u\right\|_2^2\|u\|_2^{r-2}}\right\}.
\end{aligned}
\]
By \eqref{2411192049}, we see that under the condition \eqref{2411192127}, the supremum of $Q$ is finite.

\smallskip
Remark that here we could not deduce what conditions can guarantee $R=\infty$ when $\alpha_1\neq \alpha_2$, which is an open problem in this paper.

\smallskip

\begin{proof}[\bf Proof~of~Theorem~\ref{T2025011217} completed] The above analysis and Proposition \ref{2501171415} can conclude our theorem.
\end{proof}

\section{Mass-supercritical case}\label{2501190022}
\noindent

In this section, we consider the existence of normalized solution for problem \eqref{23082801}-\eqref{23101001} with $r_{1}+r_{2}\in(\bar{r},2^{**})$. We will work in the space $H_r^2(\mathbb{R}^N)\times H_r^2(\mathbb{R}^N)$ as in \eqref{2507102258}. For the sake of the compact embedding $H_r^{2}(\mathbb{R}^N)\subset L^p(\mathbb{R}^N)$ with $2<p<2^{**}$, we will assume $N\ge 2$. Furthermore, we denote
\[
S_\rho^r=\left\{(u, v)\in S_\rho: u,\, v \mbox{ are radially symmetric}\right\}.
\]

From \eqref{23092806}, we have
\begin{equation*}
\begin{aligned}
I(u,v)\geq& \frac{1}{2}\left(\|\Delta u\|^{2}_{2}+\|\Delta v\|^{2}_{2}\right)
-\frac{1}{2}\mathcal{D}_{2}\rho\left(\|\Delta u\|^{2}_{2}+\|\Delta v\|^{2}_{2}\right)^{\frac{1}{2}}
-\mathcal{D}_{1}\beta\rho^{r\left(1-\gamma_{r}\right)}\left(\|\Delta u\|^{2}_{2}
+\|\Delta v\|^{2}_{2}\right)^{\frac{r \gamma_{r}}{2}}  \\
=:& h\left((\|\Delta u\|^{2}_{2}+\|\Delta v\|^{2}_{2})^{\frac{1}{2}}\right),
\end{aligned}
\end{equation*}
where  the function $h: \mathbb{R}^{+} \rightarrow \mathbb{R}$ is given by
\begin{equation*}
\begin{aligned}
h(t):=\frac{1}{2} t^{2}
-\frac{1}{2}\mathcal{D}_{2}\rho t
-\mathcal{D}_{1}\beta \rho^{r\left(1-\gamma_{r}\right)} t^{r\gamma_{r}}.
\end{aligned}
\end{equation*}
To understand the geometry structure of the functional $I|_{S_{\rho}}$, we shall analyze the function $h$.
\begin{lemma}\label{L23091101}
Let $\alpha_{1}$, $\alpha_{2}$, $\beta>0$, $r\in(\bar{r},2^{**})$ and $\beta \rho^{r-2}<c^{*}(N,r)$. Then there exist $0<R_{0}<R_{1}$, such that $ h\left(R_{0}\right)=h\left(R_{1}\right)=0$ and $h(t)>0$ if and only if $t \in\left(R_{0}, R_{1}\right)$.
\begin{proof}
Setting
\begin{equation*}
\begin{aligned}
\phi(t):=\frac{1}{2}t-\mathcal{D}_{1}\beta \rho^{r\left(1-\gamma_{r}\right)} t^{r\gamma_{r}-1},
\end{aligned}
\end{equation*}
we have that for $t>0$, $h(t)>0$ if and only if $\phi(t)>\frac{1}{2}\mathcal{D}_{2}\rho$. Since $r\gamma_{r}>2$, we see that $\phi(t)$ has a unique critical point $\bar{t}\in(0,+\infty)$, and it is the global maximum point, where
\begin{equation*}\label{23092203}
\begin{aligned}
\bar{t}=\left(\frac{1}{2\mathcal{D}_{1}\beta\rho^{r\left(1-\gamma_{r}\right)} (r\gamma_{r}-1)}\right)^{\frac{1}{r\gamma_{r}-2}}.
\end{aligned}
\end{equation*}
So, $\max\limits_{t>0}\phi(t)=\phi(\bar{t})=\frac{r\gamma_{r}-2}{2 (r\gamma_{r}-1)}\bar{t}$. Clearly, $\phi(+\infty)=h(+\infty)=-\infty$, $\phi(0^+)=0^{+}$ and $h(0^+)=0^{-}$, thus
\begin{equation*}
\begin{aligned}
h(\bar{t})>0~\Leftrightarrow ~ \phi(\bar{t})>\frac{1}{2}\mathcal{D}_{2}\rho ~\Leftrightarrow ~\beta \rho^{r-2}<c^{*}(N,r).
\end{aligned}
\end{equation*}
In this case, there exist $0<R_{0}<R_{1}$ such that $h>0$  if and only if $t\in(R_0,R_1)$.
\end{proof}
\end{lemma}
Under the assumption $\beta \rho^{r-2}<c^{*}(N,r)$, we denote the set
\begin{equation}\label{2501131441}
A_{\rho}:=\{(u, v)\in S^r_\rho: \|\Delta u\|^{2}_{2}+\|\Delta v\|^{2}_{2}<R^2_0\},
\end{equation}
and it is clear that there exists some $(u_0, v_0)\in A_\rho$ such that $I(u_0, v_0)<0$. Furthermore, we can easily find another point $(u_1, v_1)\in S^r_\rho\backslash A_\rho$ such that $I(u_1, v_1)<0$. Now, we can define mountain pass set as
\[
\Gamma_\rho:=\{\gamma\in C([0, 1], S^r_\rho): \gamma(0)=(u_0, v_0),\; \gamma(1)=(u_1, v_1)\},
\]
and mountain pass level
\begin{equation}\label{2501131104}
\widetilde{M}(\rho):=\underset{\gamma\in \Gamma_\rho}{\inf} \underset{t\in [0, 1]}{\max}I(\gamma(t)).
\end{equation}
It can be seen from Lemma \ref{L23091101} that $\widetilde{M}(\rho)>0$. In the sequel, we will show that the above level  $\widetilde{M}(\rho)$ is a critical value of $I|_{S_\rho^r}$. We denote the Poho\v{z}aev functional
\begin{equation*}
\begin{aligned}
P(u,v):=2\int_{\mathbb{R}^{N}}(|\Delta u|^{2}+|\Delta v|^{2}) \ud x
-\alpha_{1}\int_{\mathbb{R}^{N}}|\nabla u|^{2} \ud x
-\alpha_{2}\int_{\mathbb{R}^{N}}|\nabla v|^{2} \ud x
-2\beta r\gamma_{r}\int_{\mathbb{R}^{N}}|u|^{r_{1}}|v|^{r_{2}} \ud x.
\end{aligned}
\end{equation*}

\smallskip
By standard procedures as Jeanjean's method in \cite{Jeanjean97} and using the augmented functional
\begin{equation*}\label{2501131409}
\begin{aligned}
\Psi({(u,v)},s):=&I(s\ast(u,v))\\
=&\frac{e^{4s}}{2}\int_{\mathbb{R}^{N}}(|\Delta u|^{2}+|\Delta v|^{2}) \ud x
-\frac{\alpha_{1}e^{2s}}{2}\int_{\mathbb{R}^{N}}|\nabla u|^{2} \ud x
-\frac{\alpha_{2}e^{2s}}{2}\int_{\mathbb{R}^{N}}|\nabla v|^{2} \ud x\\
&-\beta e^{2r\gamma_{r}s}\int_{\mathbb{R}^{N}}|u|^{r_{1}}|v|^{r_{2}} \ud x,\quad \forall (u,v)\in H_r^{2}(\mathbb{R}^{N})\times H_r^{2}(\mathbb{R}^{N}),\;\; s\in \mathbb{R},
\end{aligned}
\end{equation*}
where $s\ast(u,v):=(s\ast u, s\ast v)$ and $(s\ast u)(x):=e^{\frac{N}{2}s} u\left(e^{s} x\right)$,
we can obtain a Palais-Smale sequence approaching Poho\v{z}aev manifold for the mountain pass level $\widetilde{M}(\rho)$, here we omit its proof.

\begin{lemma}\label{L24062401}
Let $\alpha_{1}$, $\alpha_{2}$, $\beta>0$, $r\in(\bar{r},2^{**})$ and $\beta \rho^{r-2}<c^{*}(N,r)$. Then there exists a Palais-Smale sequence $\{\left(u_{n},v_{n}\right)\}$ for $I|_{S^r_{\rho}}$ at $\widetilde{M}(\rho)$, which satisfies $P\left(u_{n},v_{n}\right) \rightarrow 0$.
\end{lemma}

The above lemma concludes the existence of a Palais-Smale sequence approaching Poho\v{z}aev manifold. We next analyze the compactness of such sequence.

\begin{lemma}\label{L23101001}
Assume that $\alpha_{1}$, $\alpha_{2}$, $\beta>0$, $r\in(\bar{r},2^{**})$, $N\ge 2$ and $\beta \rho^{r-2}<\min\{c^{*}(N,r),c_{*}(N,r)\}$. Let $(u_{n},v_{n})\subset S^{r}_{\rho}$ be a Palais-Smale sequence of $I|_{S_\rho^r}$ at some level $c>0$ with $P(u_{n},v_{n}) \rightarrow 0$. Then up to a subsequence, $(u_{n},v_{n})\rightarrow(u,v)$ in $H_r^{2}(\mathbb{R}^{N})\times H_r^{2}(\mathbb{R}^{N})$, and $(u,v)$ is a radial solution of problem \eqref{23082801}-\eqref{23101001} for some $\lambda>\frac{\max\{\alpha^{2}_1,\alpha^{2}_2\}}{4}$.
\begin{proof}
From $P(u_{n},v_{n}) \rightarrow 0$,  \eqref{24030301} and $r\in(\bar{r},2^{**})$, we have
\begin{align*}
c=& I(u_{n},v_{n})+o_n(1)\\
=&\frac{1}{2}(\|\Delta u_{n}\|^{2}_{2}+\|\Delta v_{n}\|^{2}_{2})
-\frac{\alpha_{1}}{2}\|\nabla u_{n}\|^{2}_2
-\frac{\alpha_{2}}{2}\|\nabla v_{n}\|^{2}_2
-\beta\int_{\mathbb{R}^{N}}|u_{n}|^{r_{1}}|v_{n}|^{r_{2}} \ud x +o_n(1)\\
=&\left(\frac{1}{2}-\frac{1}{r\gamma_{r}}\right)(\|\Delta u_{n}\|^{2}_{2}+\|\Delta v_{n}\|^{2}_{2})
-\left(\frac{1}{2}-\frac{1}{2r\gamma_{r}}\right)(\alpha_{1}\|\nabla u_{n}\|^{2}_2
+\alpha_{2}\|\nabla v_{n}\|^{2}_2)+o_n(1)\\
\geq&\left(\frac{1}{2}-\frac{1}{r\gamma_{r}}\right)(\|\Delta u_{n}\|^{2}_{2}+\|\Delta v_{n}\|^{2}_{2})
-\left(\frac{1}{2}-\frac{1}{2r\gamma_{r}}\right)\mathcal{D}_{2}\rho (\|\Delta u_{n}\|^{2}_{2}
+\|\Delta v_{n}\|^{2}_{2})^{\frac{1}{2}}+o_n(1).
\end{align*}
Combined with $(u_{n},v_{n})\subset S^{r}_{\rho}$, we can deduce that $\{(u_{n},v_{n})\}$ is bounded in $H^{2}_{r}(\mathbb{R}^{N})\times H^{2}_{r}(\mathbb{R}^{N})$.

\smallskip
By the compactness of the embedding $H^{2}_{r}(\mathbb{R}^{N}) \hookrightarrow L^{p}\left(\mathbb{R}^{N}\right)$ for $2<p<2^{**}$, there exists a $(u,v)\in H_r^{2}(\mathbb{R}^{N})\times H_r^{2}(\mathbb{R}^{N})$ such that up to a subsequence,
\begin{equation}\label{2501171703}
\begin{aligned}
\left(u_{n},v_{n}\right) \rightharpoonup& (u,v) \quad
\text { in } H^{2}_r(\mathbb{R}^{N})\times H^{2}_r(\mathbb{R}^{N}), \\
\left(u_{n},v_{n}\right) \rightarrow& (u,v) \quad
\text { in } L^{p}\left(\mathbb{R}^{N}\right)\times L^{p}\left(\mathbb{R}^{N}\right),
 \text { for } 2<p<2^{**},\\
\left(u_{n},v_{n}\right) \rightarrow& (u,v) \quad
\text { a.e.~in } \mathbb{R}^{N}\times \mathbb{R}^{N}.
\end{aligned}
\end{equation}
Since $I'|_{S^{r}_{\rho}}(u_{n},v_{n}) \rightarrow 0$, by the Lagrange multipliers rule, we get that there exists a sequence $\{\lambda_{n}\}\subset\mathbb{R}$ such that
\begin{align}\label{23101002}
I'(u_n, v_n)+\lambda_n(u_n, v_n)\to 0\quad \mbox{ in } \left(H^2_r(\mathbb{R}^N)\times H^2_r(\mathbb{R}^N)\right)^*.
\end{align}
Multiplying \eqref{23101002} by $(u_{n},v_{n})$, we have
\begin{equation}\label{23101003}
\begin{aligned}
\lambda_{n}\rho^{2}
=&-(\|\Delta u_n\|^{2}_2+\|\Delta v_n\|^{2}_2)
+\alpha_{1}\|\nabla u_n\|^{2}_2+\alpha_{2}\|\nabla v_n\|^{2}_2
+\beta r\int_{\mathbb{R}^{N}}|u_n|^{r_{1}}|v_n|^{r_{2}} \ud x+o_n(1).\\
\end{aligned}
\end{equation}
Thus by \eqref{23091001}, we deduce that $\{\lambda_{n}\}$ is bounded, and hence up to a subsequence $\lambda_{n}\rightarrow \lambda$ for some $\lambda\in\mathbb{R}$.

\medskip
Next, we claim that $u\neq 0$ and $v\neq0$. If $u=0$ or $v=0$, we deduce from \eqref{2501171703} that
\begin{equation}\label{2501171708}
\begin{aligned}
\lim_{n\rightarrow\infty}\int_{\mathbb{R}^{N}}|u_n|^{r_{1}}|v_n|^{r_{2}} \ud x=0.
\end{aligned}
\end{equation}
With \eqref{2501171708} and using $P(u_{n},v_{n}) \rightarrow 0$ again,
\begin{equation*}
\begin{aligned}
0<c=&\lim_{n\rightarrow\infty} I(u_{n},v_{n})
=\lim_{n\rightarrow\infty}\left( \frac{1}{2}(\|\Delta u_{n}\|^{2}_{2}+\|\Delta v_{n}\|^{2}_{2})
-\frac{\alpha_{1}}{2}\|\nabla u_{n}\|^{2}_2
-\frac{\alpha_{2}}{2}\|\nabla v_{n}\|^{2}_2\right)\\
=&-\frac{1}{2}\lim_{n\rightarrow\infty}\left(\|\Delta u_{n}\|^{2}_{2}+\|\Delta v_{n}\|^{2}_{2}\right)\leq 0,
\end{aligned}
\end{equation*}
which is impossible. Thus, $u\neq 0$ and $v\neq 0$.

\medskip
{\it Claim.}  $\liminf\limits_{n\rightarrow\infty}
\left(\|\Delta u_{n}\|^{2}_{2}+\|\Delta v_{n}\|^{2}_{2}\right)$ is positive and $\lambda>\frac{\max\{\alpha^{2}_1,\alpha^{2}_2\}}{4}$.

\smallskip
According to $P(u_{n},v_{n}) \rightarrow 0$, we deduce that
\begin{equation*}
\begin{aligned}
0<c= I(u_{n},v_{n})+o_n(1)
=-\frac{1}{2}(\|\Delta u_{n}\|^{2}_{2}+\|\Delta v_{n}\|^{2}_{2})
+\beta (r\gamma_{r}-1)\int_{\mathbb{R}^{N}}|u_{n}|^{r_{1}}|v_{n}|^{r_{2}} \ud x+o_n(1).
\end{aligned}
\end{equation*}
This combined with \eqref{23091001} shows that
\begin{equation}\label{23101802}
\begin{aligned}
 \frac{1}{2}(\|\Delta u_{n}\|^{2}_{2}+\|\Delta v_{n}\|^{2}_{2})
\leq&\beta (r\gamma_{r}-1)\int_{\mathbb{R}^{N}}|u_{n}|^{r_{1}}|v_{n}|^{r_{2}} \ud x+o_n(1)\\
\leq&\beta (r\gamma_{r}-1)\mathcal{D}_{1}\rho^{r(1-\gamma_r)}\left(\|\Delta u_{n}\|^{2}_{2}
+\|\Delta v_{n}\|^{2}_{2}\right)^{\frac{r \gamma_{r}}{2}}+o_n(1).
\end{aligned}
\end{equation}
From $r\in(\bar{r},2^{**})$, we obtain
\begin{equation}\label{2501171755}
\begin{aligned}
\liminf_{n\rightarrow\infty}(\|\Delta u_{n}\|^{2}_{2}+\|\Delta v_{n}\|^{2}_{2})
\geq\left(\frac{1}{2\beta (r\gamma_{r}-1)\mathcal{D}_{1}\rho^{r(1-\gamma_r)}}\right)^{\frac{2}{r \gamma_{r}-2}}.
\end{aligned}
\end{equation}

Next, we prove that $\lambda>\frac{\max\{\alpha^{2}_1,\alpha^{2}_2\}}{4}$. Firstly we get from $P\left(u_{n},v_{n}\right) \rightarrow 0$ that
\begin{equation}\label{24101402}
\begin{aligned}
\lim_{n\rightarrow\infty}\int_{\mathbb{R}^{N}}|u_{n}|^{r_{1}}|v_{n}|^{r_{2}} \ud x
\leq \frac{1}{\beta r\gamma_{r}} \liminf_{n\rightarrow\infty}(\|\Delta u_{n}\|^{2}_{2}+\|\Delta v_{n}\|^{2}_{2}).
\end{aligned}
\end{equation}
In the sequel, we prove the claim by two cases.

\smallskip
{\it Case. $\frac{1}{2}<\gamma_{r}<1$.} Using $P\left(u_{n},v_{n}\right) \rightarrow 0$, \eqref{23101003}, \eqref{2501171755} and \eqref{24101402}, we have
\begin{align*}
\lambda\rho^{2}
=&\liminf_{n\rightarrow\infty}(\|\Delta u_{n}\|^{2}_{2}+\|\Delta v_{n}\|^{2}_{2})
+\beta r(1-2\gamma_{r})\lim_{n\rightarrow\infty}\int_{\mathbb{R}^{N}}|u_n|^{r_{1}}|v_n|^{r_{2}} \ud x\\
\geq& \frac{1-\gamma_{r}}{\gamma_{r}}\liminf_{n\rightarrow\infty}(\|\Delta u_{n}\|^{2}_{2}+\|\Delta v_{n}\|^{2}_{2})\\
\geq& \frac{1-\gamma_{r}}{\gamma_{r}} \left(\frac{1}{2\beta (r\gamma_{r}-1)\mathcal{D}_{1}\rho^{r(1-\gamma_r)}}\right)^{\frac{2}{r \gamma_{r}-2}}.
\end{align*}
Since $\beta \rho^{r-2}<c_{*}(N,r)$, thus $\lambda>\frac{\max\{\alpha^{2}_1,\alpha^{2}_2\}}{4}$.

\medskip
{\it Case. $0<\gamma_{r}\leq\frac{1}{2}$.} By $P\left(u_{n},v_{n}\right) \rightarrow 0$, \eqref{23101003}, \eqref{23101802} and \eqref{2501171755}, we have
\begin{align*}
\lambda\rho^{2}
=&\liminf_{n\rightarrow\infty}(\|\Delta u_{n}\|^{2}_{2}+\|\Delta v_{n}\|^{2}_{2})
+\beta r(1-2\gamma_{r})\lim_{n\rightarrow\infty}\int_{\mathbb{R}^{N}}|u_n|^{r_{1}}|v_n|^{r_{2}} \ud x\\
\geq& \frac{r-2}{2(r\gamma_{r}-1)}\liminf_{n\rightarrow\infty}(\|\Delta u_{n}\|^{2}_{2}+\|\Delta v_{n}\|^{2}_{2})\\
\geq& \frac{r-2}{2(r\gamma_{r}-1)}\left(\frac{1}{2\beta (r\gamma_{r}-1)\mathcal{D}_{1}\rho^{r(1-\gamma_r)}}\right)^{\frac{2}{r \gamma_{r}-2}}.
\end{align*}
Since $\beta \rho^{r-2}<c_{*}(N,r)$, thus $\lambda>\frac{\max\{\alpha^{2}_1,\alpha^{2}_2\}}{4}$. We complete the proof of the claim.

\medskip
Since $u_{n}\rightharpoonup u \neq 0$ and $v_{n}\rightharpoonup v\neq 0$ weakly in $H_r^{2}(\mathbb{R}^{N})\times H_r^{2}(\mathbb{R}^{N})$ and by \eqref{23101002}, we have
\begin{equation}\label{24101403}
\begin{aligned}
I'(u,v) +\lambda(u, v)=0 \quad \mbox{in }\; \left(H_r^{2}(\mathbb{R}^{N})\times H_r^{2}(\mathbb{R}^{N})\right)^*.
\end{aligned}
\end{equation}
Let $(u_{n}-u, v_{n}-v)$ multiply \eqref{23101002} and \eqref{24101403}. By $\lambda_{n}\rightarrow \lambda$, we obtain
\begin{equation*}
\begin{aligned}
(I'(u_{n},v_{n})-I'(u,v))[u_{n}-u,v_{n}-v]+\lambda\int_{\mathbb{R}^{N}} (|u_{n}-u|^2+|v_{n}-v|^2) \ud x=o_{n}(1).
\end{aligned}
\end{equation*}
Since $(u_{n},v_{n}) \rightarrow(u,v)$ in $L^{p}(\mathbb{R}^{N}) \times L^{p}(\mathbb{R}^{N})$ for $2<p<2^{**}$ and using \eqref{24030301}, we deduce that
\begin{equation}\label{24101404}
\begin{aligned}
&\|\Delta(u_{n}-u)\|^{2}_{2}+\|\Delta(u_{n}-u)\|^{2}_{2}
+\lambda(\|u_{n}-u\|^{2}_{2}+\|v_{n}-v\|^{2}_{2})\\
=&\alpha_1\|\nabla(u_{n}-u)\|^{2}_{2}
+\alpha_2\|\nabla(v_{n}-v)\|^{2}_{2}+o_{n}(1)\\
\leq&\max\{\alpha_1,\alpha_2\}\left(\|\Delta(u_{n}-u)\|^{2}_{2}
+\|\Delta(v_{n}-v)\|^{2}_{2}\right)^{\frac{1}{2}}
(\|u_{n}-u\|^{2}_{2}+\|v_{n}-v\|^{2}_{2})^{\frac{1}{2}}+o_{n}(1).
\end{aligned}
\end{equation}

Assume that $\|\Delta(u_{n}-u)\|^{2}_{2}
+\|\Delta(v_{n}-v)\|^{2}_{2}\geq \delta$ and $\|u_{n}-u\|^{2}_{2}+\|v_{n}-v\|^{2}_{2}\geq \delta$ for some $\delta>0$. By \eqref{24101404}, we have
\begin{equation}\label{24033103}
\begin{aligned}
2 \sqrt{\lambda}\leq& \frac{\left(\|\Delta(u_{n}-u)\|^{2}_{2}
+\|\Delta(v_{n}-v)\|^{2}_{2}\right)^{\frac{1}{2}}}
{(\|u_{n}-u\|^{2}_{2}+\|v_{n}-v\|^{2}_{2})^{\frac{1}{2}}}
+\lambda \frac
{(\|u_{n}-u\|^{2}_{2}+\|v_{n}-v\|^{2}_{2})^{\frac{1}{2}}}
{\left(\|\Delta(u_{n}-u)\|^{2}_{2}+\|\Delta(v_{n}-v)\|^{2}_{2}\right)^{\frac{1}{2}}}\\
\leq& \max\{\alpha_1,\alpha_2\}+o_{n}(1),
\end{aligned}
\end{equation}
which contradicts $\lambda>\frac{\max\{\alpha^{2}_1,\alpha^{2}_2\}}{4}$. So $\|\Delta(u_{n}-u)\|^{2}_{2}
+\|\Delta(v_{n}-v)\|^{2}_{2}\rightarrow 0$ or $\|u_{n}-u\|^{2}_{2}+\|v_{n}-v\|^{2}_{2}\rightarrow 0$. In both cases, using again \eqref{24101404} and $\lambda>0$, we can prove that
\begin{equation*}
\begin{aligned}
\|\Delta(u_{n}-u)\|^{2}_{2},~
\|\Delta(v_{n}-v)\|^{2}_{2},~
\|u_{n}-u\|^{2}_{2},~
\|v_{n}-v\|^{2}_{2} \rightarrow 0 \quad \text { as } n \rightarrow \infty.
\end{aligned}
\end{equation*}
Thus, $(u_{n},v_{n}) \rightarrow (u,v)$ strongly in $H_r^{2}(\mathbb{R}^{N})\times H_r^{2}(\mathbb{R}^{N})$.
\end{proof}
\end{lemma}

\begin{proof}[\bf Proof of Theorem \ref{T230828010}  completed] Under $\beta \rho^{r-2}<c^{*}(N,r)$, let $\widetilde{M}(\rho)$ be given in \eqref{2501131104}, which is positive. By Lemma \ref{L24062401}, we obtain a Palais-Smale sequence $(u_n, v_n)$ of $I|_{S^r_{\rho}}$ at the level $\widetilde{M}(\rho)$, which satisfies $P\left(u_{n},v_{n}\right) \rightarrow 0$. Using Lemma \ref{L23101001}, if in addition $\beta \rho^{r-2}<c_{*}(N,r)$, we obtain a radial solution $(u, v)$ for some $\lambda>\frac{\max\{\alpha^{2}_1,\alpha^{2}_2\}}{4}$.
\end{proof}

%
%
%
%

\noindent {\bf \large Acknowledgements.} Z. Jin was supported by Fundamental Research Program of Shanxi Province of China(No. 202303021212160) and the China Scholarship Council(No. 202508140048), G. Wang was supported by Fundamental Research Program of Shanxi Province of China (No. 202403021221163). Last we thank the anonymous reviewer for the careful reading of our manuscript and the insightful comments and suggestions.

\end{document}